\documentclass[journal,twoside]{IEEEtran}

\usepackage{cite}
\usepackage[cmex10]{amsmath,mathtools}
\usepackage{amsfonts}
\usepackage{amssymb}
\usepackage{textcomp}
\usepackage[caption=false,font=footnotesize]{subfig}
\usepackage{fixltx2e}
\usepackage{url}
\usepackage{comment}
\usepackage{siunitx}
\usepackage{xspace}
\usepackage{accents}
\usepackage{ifthen}
\usepackage{bm}
\MakeRobust{\underaccent}

\newtheorem{theorem}{Theorem}
\newtheorem{lemma}{Lemma}

\newtheorem{corollary}{Corollary}
\newtheorem{definition}{Definition}

\newtheorem{assumption}{Assumption}

\newtheorem{remarknonum}{Remark}

\DeclareMathOperator{\real}{Re}
\DeclareMathOperator{\imag}{Im}
\DeclareMathOperator{\trace}{tr}
\DeclareMathOperator{\rank}{rank}

\DeclareMathOperator{\nullspace}{null}
\DeclareMathOperator{\dimension}{dim}
\DeclareMathOperator*{\minimize}{minimize}
\DeclareMathOperator*{\maximize}{maximize}
\DeclareMathOperator*{\subjectto}{subject~to}

\DeclareMathOperator{\epi}{epi}

\let\vec\relax\DeclareMathOperator{\vec}{vec}

\DeclareMathOperator{\interior}{int}

\newcommand{\subdiff}{\partial}
\DeclareMathOperator{\psd}{\hskip1pt\succeq\hskip1pt}
\DeclareMathOperator{\pd}{\hskip1pt\succ\hskip1pt}

\newcommand{\abs}[1]{{\lvert#1\rvert}}

\newcommand{\vm}[1]{{\bm{#1}}}
\newcommand{\mc}[1]{{\mathcal{#1}}}

\newcommand{\conj}{^\ast}
\DeclareMathOperator{\transpose}{T}
\DeclareMathOperator{\hermitian}{H}
\newcommand{\tran}{^{\transpose}}
\newcommand{\herm}{^{\hermitian}}
\newcommand{\nR}{\mathbb{R}}
\newcommand{\nRnn}{\mathbb{R}_{+}}

\newcommand{\nN}{\mathbb{N}}
\newcommand{\nC}{\mathbb{C}}
\newcommand{\nS}{\mathbb{S}}
\newcommand{\nSp}{\tilde{\mathbb{S}}}
\newcommand{\aOut}[1]{\expandafter\hat#1}
\newcommand{\aIn}[1]{\expandafter\check#1}
\newcommand{\aBar}[1]{\expandafter\bar#1}
\newcommand{\aUBar}[1]{\expandafter\underaccent{\bar}#1}
\newcommand{\aShunt}[1]{\expandafter\tilde#1}
\newcommand{\aArrow}[1]{\expandafter\vec#1}
\newcommand{\aDot}[1]{\expandafter\dot#1}
\newcommand{\aLoss}[1]{\expandafter\tilde#1}
\newcommand{\aCost}[1]{\expandafter\bar#1}
\newcommand{\aDc}[1]{\expandafter\mathring#1}

\newcommand{\aGen}[1]{#1^{\mathrm{G}}}
\newcommand{\aLoad}[1]{#1^{\mathrm{L}}}
\newcommand{\aLB}[1]{\expandafter\underaccent{\bar}#1}
\newcommand{\aUB}[1]{\expandafter\bar#1}

\newcommand{\aOutUB}[1]{\bar{\hat{#1}}}

\newcommand{\aInUB}[1]{\bar{\check{#1}}}

\newcommand{\aOpt}[1]{#1^{\star}}

\newcommand{\aArc}[1]{\expandafter\wideparen#1}

\newcommand{\sV}{{\mc{V}}}
\newcommand{\sE}{{\mc{E}}}
\newcommand{\sD}{{\mc{D}}}
\newcommand{\sB}{{\mc{B}}}

\newcommand{\sN}{{\mc{N}}}

\newcommand{\sS}{{\mc{S}}}

\newcommand{\sH}{{\mc{H}}}

\newcommand{\sW}{{\mc{W}}}

\newcommand{\sG}{{\mc{G}}}
\newcommand{\sM}{{\mc{M}}}
\newcommand{\sA}{{\mc{A}}}
\newcommand{\sI}{{\mc{I}}}
\newcommand{\gG}{{\mc{G}}}
\newcommand{\ehat}{{\aOut{\epsilon}}}
\newcommand{\echk}{{\aIn{\epsilon}}}
\newcommand{\dhat}{{\aOut{\delta}}}
\newcommand{\dchk}{{\aIn{\delta}}}
\newcommand{\fPsi}{\varPsi}
\newcommand{\fpsi}{\psi}

\newcommand{\fB}{B}
\newcommand{\fC}{C}
\newcommand{\fXi}{\Xi}

\newcommand{\lagrange}{L}
\newcommand{\fSigmaB}{\sigma}
\newcommand{\fSigmaC}{\pi}
\newcommand{\cNV}{{N}_\sV}
\newcommand{\cNE}{{N}_\sE}
\newcommand{\cND}{{N}_\sD}
\newcommand{\cM}{M}

\newcommand{\mM}{{\vm{M}}}
\newcommand{\mV}{{\vm{V}}}
\newcommand{\mI}{{\vm{I}}}
\newcommand{\mIout}{{\hat{\skew{-3}\bm{I}}\hskip-2.25pt}}
\newcommand{\mIin}{{\check{\skew{-3}\bm{I}}\hskip-2.25pt}}
\newcommand{\mS}{{\vm{S}}}
\newcommand{\mP}{{\vm{P}}}
\newcommand{\mQ}{{\vm{Q}}}
\newcommand{\mA}{{\vm{A}}}

\newcommand{\mL}{{\vm{L}}}

\newcommand{\mAub}{{\bar{\skew{-6}\bm{A}}\hskip-3.5pt}}
\newcommand{\mB}{{\vm{B}}}
\newcommand{\mC}{{\vm{C}}}

\newcommand{\mLambda}{{\vm{\varLambda}}}

\newcommand{\mZero}{{\vm{0}}}
\newcommand{\mVp}{{\tilde{\vm{V}}}}
\newcommand{\vv}{{\vm{v}}}

\newcommand{\ve}{{\vm{e}}}
\newcommand{\vx}{{\vm{x}}}

\newcommand{\vp}{{\vm{p}}}
\newcommand{\vh}{{\vm{h}}}
\newcommand{\vc}{{\vm{c}}}
\newcommand{\vd}{{\vm{d}}}
\newcommand{\vg}{{\vm{g}}}
\newcommand{\vb}{{\vm{b}}}

\newcommand{\va}{{\vm{a}}}
\newcommand{\vLambda}{{\vm{\lambda}}}
\newcommand{\vXi}{{\vm{\xi}}}
\newcommand{\vMu}{{\vm{\mu}}}
\newcommand{\iu}{{\mkern1.5mu\text{\bfseries i}\mkern1.5mu}}
\newcommand{\circledA}{\raisebox{-1.25pt}{\includegraphics{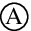}}}
\newcommand{\circledB}{\raisebox{-1.25pt}{\includegraphics{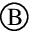}}}
\newcommand{\circledC}{\raisebox{-1.25pt}{\includegraphics{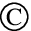}}}
\newcommand{\circledD}{\raisebox{-1.25pt}{\includegraphics{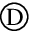}}}

\begin{document}
\title{The Hybrid Transmission Grid Architecture: Benefits in Nodal Pricing}

\author{Matthias~Hotz,~\IEEEmembership{Student Member,~IEEE,}~and~Wolfgang~Utschick,~\IEEEmembership{Senior Member,~IEEE}}

\maketitle

\begin{abstract}
Recently, we proposed a capacity expansion approach for transmission grids that combines the upgrade of transmission capacity with a transition in system structure to improve grid operation. The key to this concept is a particular hybrid AC/DC transmission grid architecture, which is obtained by uprating selected AC lines via a conversion to HVDC.~We have shown that this system structure improves optimal power flow (OPF) solvability and that it can reduce the total generation costs. In this work, we study the benefits of this hybrid architecture in the context of a deregulated electricity market. We propose an efficient and accurate nodal pricing method based on locational marginal prices (LMPs) that utilizes a second-order cone relaxation of the OPF problem. Applicability of this method requires exactness of the relaxation, which is difficult to obtain for conventional meshed AC transmission grids. We prove that the hybrid architecture ensures applicability if the LMPs do not coincide with certain pathological price profiles, which are shown to be unlikely under normal operating conditions. Using this nodal pricing method, we demonstrate that upgrading to the hybrid architecture can not only increase the effective transmission capacity but also reduce the separation of nodal markets and improve the utilization of generation.
\end{abstract}


\section{Introduction}
	\label{sec:introduction}

\IEEEPARstart{I}{n} electricity markets, \emph{nodal pricing} is an instrument to account for system constraints and losses~\cite{Kirschen2004a,Stoft2002a,Gan2013b,Jokic2007a}. The limitations on power flow in a congested grid impose restrictions on trades and electrical losses distort the balance of supply and demand, which is reflected by bus-dependent (nodal) prices for electrical power. For a market with \emph{perfect competition}, i.e., when all market participants are price takers and do not exert market power, the optimal nodal prices equal the \emph{locational marginal prices} (LMPs)~\cite{Kirschen2004a}. The marginal price of active power at a bus corresponds to the cost of serving an increment of load by the cheapest possible means of generation~\cite{Kirschen2004a}, i.e., LMPs capture the sensitivity of the minimum total generation cost to load variations. Accordingly, LMPs are tightly related to the \emph{optimal power flow} (OPF) problem, which is an optimization problem that identifies the minimum cost generation dispatch for a given load considering an AC model of the grid and system constraints. In particular, LMPs quantify the sensitivity of the optimal objective value of the OPF problem with respect to the nodal power balance~constraints.

For conventional meshed AC transmission grids, accurate LMPs are hard to obtain due to the nonconvexity of the OPF problem~\cite{Overbye2004a,Wang2007a}. As a consequence, LMPs are typically determined on the basis of a simplified system model known as ``DC power flow'' (cf. e.g.~\cite{Li2007b}), which constitutes a linearization of the AC power flow equations that considers only active power and assumes lossless lines, a flat voltage profile, and small bus voltage angle differences~\cite{Wood1996a}. The OPF formulation based on the DC power flow is known as ``DC\;OPF'' and constitutes a linear program. Due to strong duality in linear programs, the associated approximation of LMPs is given by the Lagrangian dual variables of the power balance constraints~\cite{Kirschen2004a}, which are efficiently computed, e.g., using an interior-point method~\cite{Boyd2004a}. However, due to the model mismatch these approximate LMPs are potentially inaccurate and may induce constraint-violating power flows, necessitating compensation measures~\cite{Li2007a}.

Recently, we proposed a \emph{hybrid architecture} that is established by a topology-conserving capacity expansion approach~\cite{Hotz2016a}. Therein, the capacity of certain transmission lines is uprated via a conversion to HVDC, where the lines are selected such that loops are resolved. We proved in~\cite{Hotz2016a} that the OPF problem of the resulting hybrid AC/DC grid permits an \emph{exact} semidefinite relaxation if the injection lower bounds are inactive, which enables its globally optimal solution with efficient polynomial time algorithms. Furthermore, the simulation results in~\cite{Hotz2016a} show that the hybrid architecture induces substantial \emph{flexibility} in power flow, which can enable a reduction of the total generation cost.

This work continues the study of this hybrid architecture in the context of a deregulated electricity market, where its benefit turns out to be twofold. On one hand, the hybrid architecture gives rise to a computationally efficient and accurate nodal pricing method and, on the other hand, it can reduce trading restrictions and improve grid operation. In line with our previous results, these findings encourage such a transition from conventional to structure-promoting capacity~expansion.

\subsection{Contributions and Outline}

Section~\ref{sec:model} presents the system model, which generalizes the hybrid transmission grid model in~\cite{Hotz2016a} to convex generation cost functions, arbitrary convex injection regions, and flexible loads for adequacy in a market context. Section~\ref{sec:opf} introduces the corresponding OPF problem and Section~\ref{sec:lmp} derives its relation to LMPs, which is shown to require a zero duality gap. On this basis, Section~\ref{sec:sdrpricing} discusses nodal pricing based on a semidefinite relaxation of the OPF problem, as exactness of this relaxation implies a zero duality gap. Two major issues thereof are identified, i.e., (a) computational inefficiency for large-scale grids and (b) the uncertainty about applicability due to potential inexactness of the relaxation. Issue (a) is addressed in Section~\ref{sec:socpricing}, which presents a further relaxation to a second-order cone problem in order to utilize sparsity to improve computational efficiency. This relaxation is then established as a nodal pricing method by proving that it maintains the relation to LMPs under exactness. Issue (b) is addressed in Section~\ref{sec:exactness} via a study of exactness. For conventional grids, the tendency towards exactness is difficult to characterize. In contrast, we show that the hybrid architecture gives rise to an intuitive characterization of exactness via the notion of pathological price profiles. To this end, exactness is proven to obtain as long as the LMPs do not coincide with certain pathological price profiles, which are unlikely under normal operating conditions. This result also extends our previous work on OPF in~\cite{Hotz2016a}, where exactness of a semidefinite relaxation is guaranteed for a less expressive system model and the technical requirement of excluding power injection lower bounds. In this regard, the softening of exactness under pathological price profiles can be understood as the trade-off for a more advanced system model and the consideration of power injection lower bounds.

Section~\ref{sec:simulation} presents an upgrade strategy and illustrates the application of the proposed nodal pricing method to a large-scale, real-world transmission grid. On one hand, these results illustrate that the hybrid architecture enables the efficient identification of accurate LMPs. On the other hand, they show that the hybrid architecture can not only increase the effective transmission capacity but also improve grid operation by substantially reducing grid-induced trading restrictions and facilitating a more efficient utilization of generation. Finally, Section~\ref{sec:conclusion} concludes the paper.

\subsection{Notation}

The set of natural numbers is denoted by $\nN$, the set of real numbers by $\nR$, the set of nonnegative real numbers by $\nRnn$, the set of complex numbers by $\nC$, and the set of Hermitian matrices in $\nC^{N\times N}$ by $\nS^N$. The imaginary unit is denoted by $\displaystyle\iu=\sqrt{-1}$. For $x\in\nC$, its real part is $\real(x)$, its imaginary part is $\imag(x)$, its absolute value is $\abs{x}$, and its complex conjugate is $x\conj$. For a matrix $\mA$, its transpose is $\mA\tran$, its conjugate (Hermitian) transpose is $\mA\herm$, its trace is $\trace(\mA)$, its rank is $\rank(\mA)$, its nullspace (kernel) is $\nullspace(\mA)$, its element in row $i$ and column $j$ is $[\mA]_{i,j}$, and its vectorization is $\vec(\mA)=[\va_1\tran,\ldots,\va_N\tran]\tran$, where $\va_1$ to $\va_N$ are the columns of $\mA$. For a complex-valued matrix $\mM=\mA+\iu\mB\in\nC^{M\times N}$, where $\mA,\mB\in\nR^{M\times N}$, its real part is $\real(\mM)=\mA$ and its imaginary part is $\imag(\mM)=\mB$. For two matrices $\mA,\mB\in\nS^N$, $\mA\psd\mB$ denotes that $\mA-\mB$ is positive semidefinite and $\mA\pd\mB$ that $\mA-\mB$ is positive definite. For real-valued vectors, inequalities are considered component-wise. The vector $\ve_n$ denotes the $n$th standard basis vector of appropriate dimension. For a set $\sS$, its cardinality is denoted by $\abs{\sS}$ and its interior by $\interior(\sS)$. For two sets $\sA$ and $\sB$, $\sA+\sB$ denotes the Minkowski sum and $\sA-\sB$ the Minkowski difference of $\sA$ and $\sB$. For a set $\sN\subset\nN$ and vectors or matrices $\vx_n\in\sS$, with $n\in\sN$, $\vx_{\sN}$ denotes the $\abs{\sN}$-tuple $\vx_{\sN} = (\vx_n)_{n\in\sN}$ and $\sS_{\sN}$ the $\abs{\sN}$-fold Cartesian product $\sS_{\sN}={\prod_{n\in\sN} \sS}$. For a vector space $\sW$, its dimension is $\dimension(\sW)$. For a convex function $f:\sS\rightarrow\nR$, $\epi(f)$ denotes the epigraph of $f$, and, for and $\vx\in\sS$, $\subdiff f(\vx)$ denotes the subdifferential of $f$ at $\vx$ and, if $f$ is differentiable at $\vx$, $\nabla f(\vx)$ denotes the gradient of $f$ at $\vx$.

\section{System Model}
	\label{sec:model}

This work utilizes the system model for hybrid transmission grids in~\cite{Hotz2016a}, which comprises an architectural and electrical part. The architecture is described by the directed multigraph $\gG = {\{\sV, \sE, \sD, \ehat, \echk, \dhat, \dchk\}}$, where $\sV = \{1,\ldots,\cNV\}$ is the set of buses, $\sE = \{1,\ldots,\cNE\}$ the set of \emph{AC branches} (AC lines, cables, transformers, and phase shifters), and $\sD = \{1,\ldots,\cND\}$ the set of \emph{DC branches} (HVDC lines, cables, and back-to-back converters). The functions $\ehat,\echk:\sE\rightarrow\sV$ and $\dhat,\dchk:\sD\rightarrow\sV$ map an AC and DC~branch to its source and destination bus, respectively. The electrical behavior is described by a steady-state model, where the system state comprises the \emph{bus voltage vector} $\vv\in\nC^{\cNV}$ and \emph{DC branch flow vector} $\vp\in\nR^{\cND}$. At the buses, generators with a rectangular injection region and fixed loads are considered. The corresponding constraints are implemented as upper and lower bounds on the net injection of active and reactive power into the grid, i.e.,
\begin{subequations}\label{eqn:model:crt:injbox}%
\begin{align}
	\aLB{P_n} \leq \vv\herm\mP_n\vv &+ \vh_n\tran\vp \leq \aUB{P_n},
		&& \forall n\in\sV\\
	\aLB{Q_n} \leq \vv\herm&\mQ_n\vv \leq \aUB{Q_n},
		&& \forall n\in\sV.
\end{align}
\end{subequations}
Therein, the matrices $\mP_n,\mQ_n\in\nS^{\cNV}$ are a function of the bus admittance matrix and characterize the flow on AC branches, while $\vh_n\in\nR^{\cND}$ describes the flow on DC branches, see~\cite[Sec.~II]{Hotz2016a}. The voltage magnitude at every bus is restricted to the corresponding voltage range $[\aLB{V_n},\aUB{V_n}]\subset\nRnn$ by
\begin{equation}\label{eqn:model:crt:vm}
	\qquad\!
	\aLB{V_n}^2 \leq \vv\herm\mM_n\vv \leq \aUB{V_n}^2,
	\qquad\qquad
		\forall n\in\sV,
\end{equation}
where $\mM_n = \ve_n\ve_n\tran\in\nS^{\cNV}$. Complementary to the definition in~\cite{Hotz2016a}, the lower bounds are assumed to be \emph{strictly} positive, i.e., ${\aLB{V_n}>0}$ for all $n\in\sV$. For AC branches, the formulation in~\cite{Hotz2016a} exchanges the usual apparent power flow limit (``MVA rating'') by its underlying constraints (see e.g.~\cite[Ch.~6.1.12]{Kundur1994a}) to improve expressiveness and mathematical structure. These constraints comprise upper bounds $\aOutUB{I}_k$ and $\aInUB{I}_k$ on the current magnitude at the source and destination, i.e.,
\begin{equation}\label{eqn:model:crt:cm}
	\vv\herm\mIout_k\vv \leq \aOutUB{I}_k^2,
		\qquad
	\vv\herm\mIin_k\vv \leq \aInUB{I}_k^2,
		\qquad\!\!
		\forall k\in\sE,
\end{equation}
a restriction of the relative bus voltage magnitude drop to the range $[\aLB{\nu_k},\aUB{\nu_k}]\subset[-1,\infty)$, i.e.,
\begin{equation}\label{eqn:model:crt:vd}
	\vv\herm\aLB{\mM_k}\vv \leq 0,
		\qquad
	\vv\herm\aUB{\mM_k}\vv \leq 0,
		\qquad\!\!\!
		\forall k\in\sE,
\end{equation}
as well as a limitation of the bus voltage angle difference to the range $[\aLB{\delta_k},\aUB{\delta_k}]\subset(-\pi/2,\pi/2)$, i.e.,
\begin{align}\label{eqn:model:crt:va}
	\hspace{-3.5mm}
	\vv\herm\mA_k\vv \leq 0,
		\ \
	\vv\herm\aLB{\mA_k}\vv \leq 0,
		\ \
	\vv\herm\mAub_k\vv \leq 0,
		&& \forall k\in\sE.
\end{align}
The current constraint matrices $\mIout_k,\mIin_k\in\nS^{\cNV}$, the voltage drop constraint matrices $\aLB{\mM_k},\aUB{\mM_k}\in\nS^{\cNV}$, as well as the angle difference constraint matrices $\mA_k,\aLB{\mA_k},\mAub_k\in\nS^{\cNV}$ are defined in~\cite[Sec.~III]{Hotz2016a}. For DC branches, upper and lower bounds on power flow are captured by the vector-valued constraint
\begin{equation}\label{eqn:model:crt:p}
	\aLB{\vp} \leq \vp \leq \aUB{\vp}\,.
\end{equation}

To further improve expressiveness, the model is extended to flexible loads (elastic demand) and a more elaborate characterization of generation. To this end, a \emph{generation vector} $\vg_n=[\aGen{P_n}, \aGen{Q_n}]\tran\in\nR^2$ and a \emph{load vector} $\vd_n=[\aLoad{P_n}, \aLoad{Q_n}]\tran{\in\nR^2}$ is introduced at every bus $n\in\sV$, where $\aGen{P_n}$ and $\aGen{Q_n}$ is the active and reactive power generation and $\aLoad{P_n}$ and $\aLoad{Q_n}$ is the active and reactive load. The P-Q capability of the (aggregated) generator at bus $n\in\sV$ is specified by the nonempty, compact, and convex set $\sG_n\subset\nR^2$, i.e., all $\vg_n\in\sG_n$ are valid operating points. For example, $\sG_n$ may constitute a simple box constraint or, for more precise modeling, the convex hull of the generator capability curve. Similarly, the admissible range for the (aggregated) load at bus $n\in\sV$ is specified by the nonempty, compact, and convex set $\sD_n\subset\nR^2$, i.e., all $\vd_n\in\sD_n$ are valid load configurations. For example, $\sD_n$ is a singleton set for a fixed load and non-singleton for a flexible load. At bus $n\in\sV$, this characterization of generation and load is linked to the hybrid transmission grid model via the power balance equations
\begin{subequations}\label{eqn:model:crt:powbal}
\begin{align}
	\vv\herm\mP_n\vv + \vh_n\tran\vp &= \ve_1\tran(\vg_n - \vd_n),
		&& \forall n\in\sV
		\label{eqn:model:crt:powbal:p}\\
	\vv\herm\mQ_n\vv &= \ve_2\tran(\vg_n - \vd_n),
		&& \forall n\in\sV,
		\label{eqn:model:crt:powbal:q}
\end{align}
\end{subequations}
which replace the power injection constraint~\eqref{eqn:model:crt:injbox}.

In a managed spot market, the market operator collects the bids and offers of producers and consumers, respectively, to set the nodal prices and clear the market~\cite{Kirschen2004a}. The bids and offers are in general increasing and decreasing staircase-shaped price-power curves that, by integration, translate to convex and concave piecewise linear cost and benefit functions for producers and consumers, respectively. This is considered by generalizing the linear generation cost in~\cite{Hotz2016a} to convex \emph{producer cost functions} $\fC_n:\sG_n\rightarrow\nR$, with $n\in\sV$. Additionally, concave \emph{consumer benefit functions} $\fB_n:\sD_n\rightarrow\nR$ are introduced to quantify the benefit perceived by flexible loads.

\begin{remarknonum}
The system model in~\cite{Hotz2016a} is introduced for hybrid transmission grids that feature the \emph{hybrid architecture}, which comprises a \emph{tree} topology of the AC subgrid as established by~\cite[Def.~6]{Hotz2016a}. If this definition is excluded, the model is suitable for hybrid transmission grids of arbitrary topology. In the following, this fact is utilized to discuss nodal pricing for general hybrid transmission grids, while the hybrid architecture is then considered from Section~\ref{sec:exactness} onwards.
\end{remarknonum}

\section{Optimal Power Flow}
	\label{sec:opf}

LMPs are related to the OPF problem by means of the sensitivity of its optimal objective value to the power balance constraints. To establish this relation, the OPF is formulated with respect to the maximum \emph{economic welfare}, i.e., the difference of the consumer's benefit and the producer's cost for a certain system state (cf. e.g.~\cite{Kirschen2004a}). With the system model above, the corresponding OPF problem reads as follows.
\begin{subequations}\label{eqn:opf:nonconvex:full}%
\begin{alignat}{2}
\aOpt{p} =\;&\maximize_{\hspace{2mm}\mathclap{\substack{\\
		\vv\in\nC^{\cNV},\,\vp\in\nR^{\cND}\\[0.075em]
		\vg_n\in\sG_n,\,\vd_n\in\sD_n}}}
	\qquad &&
	\sum_{n\in\sV} \fB_n(\vd_n) - \sum_{n\in\sV} \fC_n(\vg_n)\\[0.25em]
&\subjectto &&
\eqref{eqn:model:crt:vm},~\eqref{eqn:model:crt:cm},~\eqref{eqn:model:crt:vd},~\eqref{eqn:model:crt:va},~\eqref{eqn:model:crt:p},~\eqref{eqn:model:crt:powbal}\,.
\end{alignat}
\end{subequations}
This is a \emph{nonconvex} optimization problem due to the power balance equations and the indefiniteness of certain constraint matrices (cf.~\cite{Hotz2016a}). The nonconvexity does not only render the problem hard to solve, but also leads to a potentially nonzero duality gap with respect to the Lagrangian dual problem. As shown below, a nonzero duality gap invalidates the coupling between LMPs and the OPF problem, rendering them hard to identify.

In what follows, it is assumed that~\eqref{eqn:opf:nonconvex:full} is strictly feasible.
\begin{assumption}\label{ass:strictfeas}
There exists a feasible tuple $(\vv,\vp,\vg_{\sV},\vd_{\sV})$ in~\eqref{eqn:opf:nonconvex:full} for which~\eqref{eqn:model:crt:vm} to~\eqref{eqn:model:crt:p} hold with strict inequality, ${\vg_n-\vd_n}\in\interior(\sG_n-\sD_n)$, and $\sG_n$ and $\sD_n$ are polyhedral sets, for all $n\in\sV$.
\end{assumption}

Furthermore, the notation is condensed to simplify the exposition. Let $\mB_n=[\vec(\mP_n\tran),\vec(\mQ_n\tran)]\tran$, $\aBar{\mB}_n=[\vh_n,\mZero]\tran$, and consider that quadratic terms in $\vv$ permit the reformulation
\begin{equation}
	\vv\herm\mA\vv
		= \trace(\mA\vv\vv\herm)
		= \vec(\mA\tran)\tran\vec(\vv\vv\herm)\,.
\end{equation}
Therewith, the OPF problem~\eqref{eqn:opf:nonconvex:full} is expressed equivalently as
\begin{subequations}\label{eqn:opf:nonconvex:compact}
\begin{align}
\aOpt{p} =\;&\maximize_{
	\mathclap{\hspace{2mm}\substack{\\
		\vv\in\nC^{\cNV},\,\vp\in\nR^{\cND}\\[0.075em]
		\vg_n\in\sG_n,\,\vd_n\in\sD_n}}}
	\qquad
	\sum_{n\in\sV} \fB_n(\vd_n) - \sum_{n\in\sV} \fC_n(\vg_n)
	\\[0.25em]
&
\subjectto\notag\\
&\qquad
\mB_n\mkern-1mu\vec(\vv\vv\herm) + \aBar{\mB}_n\vp = \vg_n - \vd_n,
	&& \hspace{-7mm} n\in\sV
	\label{eqn:opf:nonconvex:compact:powbal}\\
&\qquad
\mC\mkern-0.5mu\vec(\vv\vv\herm) + \aBar{\mC}\vp \leq \vb\,.
	\label{eqn:opf:nonconvex:compact:crt:ieq}
\end{align}
\end{subequations}
Note that~\eqref{eqn:model:crt:powbal} is implemented by~\eqref{eqn:opf:nonconvex:compact:powbal} and~\eqref{eqn:model:crt:vm} to~\eqref{eqn:model:crt:p} by~\eqref{eqn:opf:nonconvex:compact:crt:ieq}, where $\mC=[\vec(\mC_1\tran),\ldots,\vec(\mC_{\cM}\tran)]\tran$, $\aBar{\mC}=[\vc_1,\ldots,\vc_{\cM}]\tran$, and $\vb=[b_1,\ldots,b_M]\tran$ are parametrized correspondingly to reproduce the $\cM=2\cNV+7\cNE+2\cND$ inequality constraints.

\section{Locational Marginal Prices}
	\label{sec:lmp}

This section illustrates in what manner and under which conditions LMPs emerge from the Lagrangian dual of the OPF problem, which serves as a basis for the nodal pricing method later on. To this end, the Lagrangian dual of the OPF problem~\eqref{eqn:opf:nonconvex:full} is derived on the basis of~\eqref{eqn:opf:nonconvex:compact}. Let $\sM = \{1,\ldots,\cM\}$ and let $\fPsi:\nR^2_{\sV}\times\nR^M\rightarrow\nS^{\cNV}$ and $\fpsi:\nR^2_{\sV}\times\nR^M\rightarrow\nR^{\cND}$ be defined as 
\begin{align}
\hspace{-0.65em}
	\fPsi(\vLambda_\sV,\vMu)
		&= \sum_{\mathclap{n\in\sV}} \Big( [\vLambda_n]_1\mP_n + [\vLambda_n]_2\mQ_n\Big)
			+ \sum_{\mathclap{m\in\sM}} [\vMu]_m\mC_m
		\label{eqn:psimtxdef}\\[0.75ex]
\hspace{-0.65em}
	\fpsi(\vLambda_\sV,\vMu)
		&= \sum_{\mathclap{n\in\sV}} [\vLambda_n]_1\vh_n + \sum_{\mathclap{m\in\sM}} [\vMu]_m\vc_m\,.
\end{align}
Therewith, the Lagrangian function of~\eqref{eqn:opf:nonconvex:compact} can be stated as
\begin{align}
\hspace{-0.75em}
	\lagrange
		&= \sum_{n\in\sV} \Big( \fB_n(\vd_n) - \vLambda_n\tran\vd_n \Big)
			+ \sum_{n\in\sV} \Big( \vLambda_n\tran\vg_n - \fC_n(\vg_n) \Big) \notag\\
		&\qquad
			- \vv\herm\fPsi(\vLambda_\sV,\vMu)\vv
			- \fpsi(\vLambda_\sV,\vMu)\tran\vp + \vMu\tran\vb
		\label{eqn:opf:nonconvex:lagrangian}
\end{align}
in which~\eqref{eqn:opf:nonconvex:compact:powbal} is dualized using the dual variables $\vLambda_n\in\nR^2$, with $n\in\sV$, and~\eqref{eqn:opf:nonconvex:compact:crt:ieq} using $\vMu\in\nRnn^M$, while the vectorization is retracted. With~\eqref{eqn:opf:nonconvex:lagrangian}, the \emph{dual problem} of~\eqref{eqn:opf:nonconvex:compact} and, thus, of the OPF problem~\eqref{eqn:opf:nonconvex:full} is obtained as
\begin{subequations}\label{eqn:opf:nonconvex:dual}
\begin{alignat}{2}
\aOpt{d} =\;&\minimize_{
	\mathclap{\hspace{2mm}\vLambda_n\in\nR^2,\,\vMu\in\nRnn^{\cM}}}
	\quad\;\ &&
	\vMu\tran\vb
		+ \sum_{n\in\sV} \fSigmaB_n(\vLambda_n)
		+ \sum_{n\in\sV} \fSigmaC_n(\vLambda_n)\\[0.25em]
&\subjectto &&
\fPsi(\vLambda_\sV,\vMu)\psd\mZero
	\label{eqn:opf:nonconvex:dual:crt:psimtx}\\
&&&\fpsi(\vLambda_\sV,\vMu)=\mZero\,.
\end{alignat}
\end{subequations}
Therein, $\fSigmaB_n:\nR^2\rightarrow\nR$ and $\fSigmaC_n:\nR^2\rightarrow\nR$ are given by
\begin{align}
	\fSigmaB_n(\vLambda_n)
		&= \max_{\vd_n\in\sD_n} \big\{ \fB_n(\vd_n) - \vLambda_n\tran\vd_n \big\}
		\label{eqn:surplus:consumer}\\[0.5ex]
	\fSigmaC_n(\vLambda_n)
		&= \max_{\vg_n\in\sG_n} \big\{ \vLambda_n\tran\vg_n - \fC_n(\vg_n) \big\}\,.
		\label{eqn:surplus:producer}
\end{align}
To relate~\eqref{eqn:opf:nonconvex:dual} to LMPs, let $(\aOpt{\vv},\aOpt{\vp},\aOpt{\vg_{\sV}},\aOpt{\vd_{\sV}},\aOpt{\vLambda_{\sV}},\aOpt{\vMu})$ be a primal and dual optimal solution of~\eqref{eqn:opf:nonconvex:full} and~\eqref{eqn:opf:nonconvex:dual}. If the \emph{duality gap} is zero (i.e., \emph{strong duality} holds), this constitutes a saddle point of the Lagrangian function~\cite{Bazaraa2006a,Boyd2004a},  which implies that $\aOpt{\vLambda_n}\in\subdiff\fC_n(\aOpt{\vg_n})$ and, if $\fC_n$ is differentiable at $\aOpt{\vg_n}$, that $\aOpt{\vLambda_n}=\nabla\fC_n(\aOpt{\vg_n})$.\footnote{Furthermore, it also holds that ${\aOpt{\vLambda_n}\in\subdiff\fB_n(\aOpt{\vd_n})}$ and, if $\fB_n$ is differentiable at $\aOpt{\vd_n}$, that $\aOpt{\vLambda_n}=\nabla\fB_n(\aOpt{\vd_n})$.} Consequently, the vector $\aOpt{\vLambda_n}$ indeed comprises the LMP for active and reactive power at bus $n\in\sV$. With the interpretation of $\vLambda_n$ as a price vector, \eqref{eqn:surplus:consumer} can be identified as the \emph{consumer's surplus function} and~\eqref{eqn:surplus:producer} as the \emph{producer's profit function}. This explains the suitability of $\aOpt{\vLambda_{\sV}}$ as \emph{nodal prices} in case the producers and consumers act profit- and surplus-maximizing, because they incentivize a welfare-maximizing behavior as $\aOpt{\vg_n}$ and $\aOpt{\vd_n}$ are maximizers in~\eqref{eqn:surplus:consumer} and~\eqref{eqn:surplus:producer}. On the contrary, if the duality gap is nonzero, the primal-dual optimal solution does \emph{not} constitute a saddle point of the Lagrangian and $\aOpt{\vLambda_{\sV}}$ may not match the LMPs, i.e., the coupling between the dual variables and LMPs is invalidated.

\begin{figure}[!t]
\centering
\includegraphics{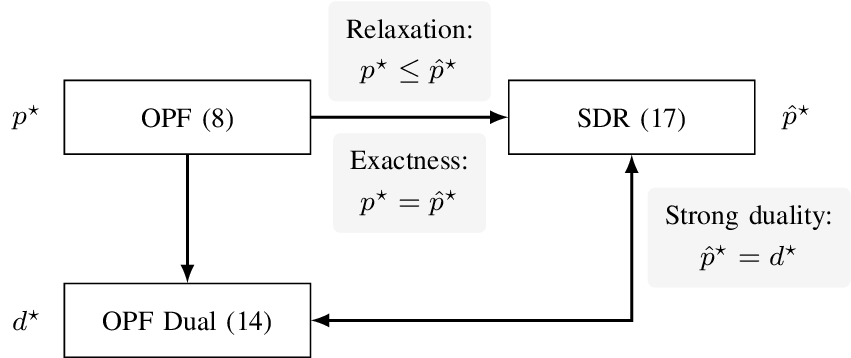}%
\caption{Interrelation of the OPF problem~\eqref{eqn:opf:nonconvex:full}, its Lagrangian dual~\eqref{eqn:opf:nonconvex:dual}, and its semidefinite relaxation~\eqref{eqn:opf:sdp}. If the relaxation is exact, it follows that $p^\star=d^\star$, i.e., the OPF problem exhibits a zero duality gap and the LMPs are given by the optimal dual variables $\aOpt{\vLambda_{\sV}}$ in~\eqref{eqn:opf:nonconvex:dual}.}%
\label{fig:relaxation:relationsdr}%
\end{figure}

\section{Nodal Pricing using Semidefinite Relaxation}
	\label{sec:sdrpricing}

A \emph{convex relaxation} of the OPF problem is thus not only motivated by computational advantages, but also by its relation to a zero duality gap. In the following, this is illustrated using the well-known semidefinite relaxation.

\subsection{Semidefinite Relaxation}
	\label{sec:sdrpricing:relaxation}

Semidefinite relaxation (SDR) is an established technique for the convex relaxation of nonconvex quadratic optimization problems~\cite{Luo2010a} and has been applied to various OPF formulations, see e.g.~\cite{Low2014a,Taylor2015a} and the references therein. Here, SDR is applied analogous to~\cite{Hotz2016a}, i.e., a Hermitian matrix $\mV\in\nS^{\cNV}$ is introduced and the quadratic expressions in $\vv$ are rewritten in terms of $\mV$. For equivalence of the formulations, $\mV$ must be positive semidefinite (psd) and have rank~$1$. In SDR, the optimization problem is rendered convex by excluding the rank constraint. Therefore, the SDR of the OPF problem~\eqref{eqn:opf:nonconvex:full} using the notation in~\eqref{eqn:opf:nonconvex:compact} reads
\begin{subequations}\label{eqn:opf:sdp}
\begin{align}
\aOpt{\hat{p}} =\;&\maximize_{
	\mathclap{\hspace{2mm}\substack{\\
		\mV\in\nS^{\cNV},\,\vp\in\nR^{\cND}\\[0.075em]
		\vg_n\in\sG_n,\,\vd_n\in\sD_n}}}
	\qquad
	\sum_{n\in\sV} \fB_n(\vd_n) - \sum_{n\in\sV} \fC_n(\vg_n)
	\\[0.25em]
&
\subjectto\notag\\
&\qquad
\mB_n\mkern-1mu\vec(\mV) + \aBar{\mB}_n\vp = \vg_n - \vd_n,
	&& \hspace{-7mm} n\in\sV\\
&\qquad
\mC\mkern-0.5mu\vec(\mV) + \aBar{\mC}\vp \leq \vb\\
&\qquad
\mV\psd\mZero\,.
	\label{eqn:opf:sdp:crt:vpsd}
\end{align}
\end{subequations}
The SDR is \emph{exact} if there exists an optimizer $(\aOpt{\mV},\aOpt{\vp},\aOpt{\vg_{\sV}},\aOpt{\vd_{\sV}})$ in~\eqref{eqn:opf:sdp} for which $\aOpt{\mV}$ has rank $1$, i.e., it facilitates the decomposition $\aOpt{\mV}={\aOpt{\vv}(\aOpt{\vv})\herm}$. Then, by construction of the SDR, it follows that $(\aOpt{\vv},\aOpt{\vp},\aOpt{\vg_{\sV}},\aOpt{\vd_{\sV}})$ is an optimizer of the OPF problem~\eqref{eqn:opf:nonconvex:full} and that $\aOpt{p}=\aOpt{\hat{p}}$.

\subsection{Nodal Pricing and Application Issues}
	\label{sec:sdrpricing:method}

A particular property of the SDR~\eqref{eqn:opf:sdp} is that its Lagrangian dual is given by~\eqref{eqn:opf:nonconvex:dual}, i.e., the OPF problem~\eqref{eqn:opf:nonconvex:full} and its SDR~\eqref{eqn:opf:sdp} share the same dual problem. Due to the convexity of the set of psd matrices, \eqref{eqn:opf:sdp}~is a \emph{convex} optimization problem and, as established by Theorem~\ref{thm:strictfeas:sdp}, Slater's constraint qualification is fulfilled, thus \emph{strong duality} holds and $\aOpt{\hat{p}}=\aOpt{d}$.
\begin{theorem}\label{thm:strictfeas:sdp}
Consider the SDR~\eqref{eqn:opf:sdp} of the OPF problem~\eqref{eqn:opf:nonconvex:full}. If Assumption~\ref{ass:strictfeas} holds, there exists a feasible tuple $(\mV,\vp,\vg_{\sV},\vd_{\sV})$ in~\eqref{eqn:opf:sdp} for which $\mV\pd\mZero$.
\end{theorem}
\begin{IEEEproof}
See Appendix~\ref{apx:thm:strictfeas:sdp}.
\end{IEEEproof}
It follows that exactness of the SDR implies a zero duality gap of the OPF problem (see also~\cite{Lavaei2012a}) and, thus, it qualifies the optimal dual variables $\aOpt{\vLambda_{\sV}}$ as LMPs, cf. Fig.~\ref{fig:relaxation:relationsdr}. Considering that interior-point methods can jointly solve the SDR~\eqref{eqn:opf:sdp} and its dual~\eqref{eqn:opf:nonconvex:dual} in polynomial time, this appears as an attractive nodal pricing method. However, there are two major issues:
\begin{itemize}
\item[(a)] \emph{Practical tractability:} In SDR, the number of optimization variables increases by $\cNV(\cNV-2)$.\footnote{The real and imaginary part are considered as separate variables and conjugate symmetry is respected. $\cNV\geq 2$.} For large-scale grids, this quadratic increase in dimensionality entails substantial difficulties in solving~\eqref{eqn:opf:sdp} and~\eqref{eqn:opf:nonconvex:dual} with an interior-point solver. This renders it highly inefficient and potentially even intractable for practical grids.
\vspace{0.25em}
\item[(b)] \emph{Applicability:} This nodal pricing method is only applicable if the semidefinite relaxation is exact. While exactness can be determined \emph{a posteriori}, the uncertainty about applicability compromises the practical value of this approach. Hence, a characterization of exactness is desired to assess if a given transmission grid tends toward an exact relaxation and, therewith, promotes applicability.
\end{itemize}

Regarding (a), a further relaxation to a second-order cone program (SOCP) was considered in the literature in order to exploit sparsity, see~\cite{Low2014a,Taylor2015a} and the references therein. This substantially reduces the computational effort for large-scale grids, but it also invalidates the interrelations in Fig.~\ref{fig:relaxation:relationsdr} and, therewith, this nodal pricing approach.

Regarding (b), there are several studies on the exactness of SDR for \emph{AC grids}. A prominent work is by Lavaei and Low~\cite{Lavaei2012a}, where a series of modified OPF problems is presented and related to exactness in order to establish an intuition about the observation of a zero duality gap. For the case of a general AC grid in~\cite[Sec.~IV-C]{Lavaei2012a}, the authors combine an extensive chain of arguments with several assumptions to convey an intuition about the range of a weighting factor on the off-diagonal block of a certain matrix (their analog of $\fPsi(\vLambda_\sV,\vMu)$), where, in case this range includes~$1$, exactness is ensured under the posed assumptions. This result is remarkable, but potentially too intricate to assess the tendency towards exactness for a specific grid. For radial grids~\cite{Low2014b} and meshed grids with loops limited to three lines~\cite{Madani2015a}, exactness has been proven under certain technical conditions, e.g., the omission of power injection lower bounds. However, these structures are not observed in transmission grids and the technical conditions cannot be maintained in practice. For \emph{hybrid} transmission grids, we recently established exactness for the hybrid architecture in~\cite{Hotz2016a}, which covers \emph{meshed} grid topologies but, still, requires the omission of power injection lower bounds. Considering that these results do \emph{not} transfer to the generalized system model in Section~\ref{sec:model}, they neither address (b) appropriately.

In the following, these two issues are investigated to arrive at an efficient nodal pricing method for large-scale transmission grids that facilitates the utilization of the \emph{hybrid architecture's} structural features to promote applicability to this system class.

\section{Nodal Pricing using SOC Relaxation}
	\label{sec:socpricing}

A major issue of nodal pricing using SDR is the quadratic increase in optimization variables, which compromises computational tractability for large-scale grids. In the following, this issue is addressed via a further relaxation to an SOCP.

\subsection{Second-Order Cone Relaxation}
	\label{sec:socpricing:relaxation}

The constraint matrices in~\eqref{eqn:model:crt:vm} to~\eqref{eqn:model:crt:powbal} are highly sparse, see~\cite[Appendix~B]{Hotz2016a}. In the following, a second-order cone (SOC) relaxation is applied to the SDR~\eqref{eqn:opf:sdp}, which enables the utilization of this sparsity to reduce the uplift in dimensionality. To this end, it is observed that a necessary (but not sufficient) condition for~\eqref{eqn:opf:sdp:crt:vpsd} is that all {$2$$\times$$2$} principal submatrices of $\mV$ are psd, cf.~\cite{Horn2013a,Kim2003a}. Thus, \eqref{eqn:opf:sdp:crt:vpsd} may be relaxed to psd constraints on {$2$$\times$$2$} principal submatrices, which can be implemented as SOC constraints, cf.~\cite{Kim2003a}. It follows from~\cite[Appendix~B]{Hotz2016a} that the sparsity pattern of the constraint matrices in~\eqref{eqn:model:crt:vm} to~\eqref{eqn:model:crt:powbal} is defined by the graph of the AC subgrid, i.e., only elements in row~$i$ and column~$j$, where $i,j\in\{\ehat(k),\echk(k)\}$ and $k\in\sE$, may be nonzero. Therefore, only the corresponding {$2$$\times$$2$} principal submatrices of $\vv\vv\herm$ in~\eqref{eqn:opf:nonconvex:compact} and, thus, of $\mV$ in~\eqref{eqn:opf:sdp} are involved in system constraints. For the SOC relaxation, this enables a reduction in dimensionality by $\cNV(\cNV - 1) - 2\cNE$ via the exclusion of irrelevant optimization variables. To formulate the SOC relaxation of~\eqref{eqn:opf:sdp}, let the set $\nSp^{\cNV}$ of Hermitian partial matrices on the AC subgraph of $\gG$ be
\begin{equation}
	\nSp^{\cNV} = \big\{ \mV\in\nS^{\cNV} : [\mV]_{i,j}=0,\ \forall(i,j)\in\{\sV\times\sV\}\setminus\sI \,\big\}
\end{equation}
where $\sI = {\{ (i,j)\in\sV\times\sV : i,j\in\{\ehat(k),\echk(k)\}, k\in\sE \}}$. Furthermore, let $\fXi_k:\nSp^{\cNV}\rightarrow\nS^2$ select the {$2$$\times$$2$} principal submatrix associated with AC branch $k\in\sE$, i.e.,
\begin{align}
	\fXi_k(\mVp) = \mS_k\tran\mVp\mS_k
\end{align}
where $\mS_k=[\ve_{\ehat(k)}, \ve_{\echk(k)}]\in\{0,1\}^{\cNV\times 2}$. Therewith, the SOC relaxation of~\eqref{eqn:opf:sdp} can be stated as
\begin{subequations}\label{eqn:opf:socp}
\begin{align}
\aOpt{\tilde{p}} =\;&\maximize_{
	\mathclap{\hspace{2mm}\substack{\\
		\mVp\in\nSp^{\cNV},\,\vp\in\nR^{\cND}\\[0.075em]
		\vg_n\in\sG_n,\,\vd_n\in\sD_n}}}
	\qquad
	\sum_{n\in\sV} \fB_n(\vd_n) - \sum_{n\in\sV} \fC_n(\vg_n)
	\\[0.25em]
&
\subjectto\notag\\
&\qquad
\mB_n\mkern-1mu\vec(\mVp) + \aBar{\mB}_n\vp = \vg_n - \vd_n,
	&& \hspace{-7mm} n\in\sV
	\label{eqn:opf:socp:crt:powbal}\\
&\qquad
\mC\mkern-0.5mu\vec(\mVp) + \aBar{\mC}\vp \leq \vb
	\label{eqn:opf:socp:crt:ineq}\\
&\qquad
\fXi_k(\mVp)\psd\mZero,
	&& \hspace{-7mm} k\in\sE\,.
	\label{eqn:opf:socp:crt:psd}
\end{align}
\end{subequations}
In~\eqref{eqn:opf:socp}, the uplift in dimensionality is reduced to a minor linear increase by $2\cNE-\cNV$ compared to the OPF problem~\eqref{eqn:opf:nonconvex:full}. Due to the further relaxation, a solution of the OPF problem can only be recovered under the following condition.
\begin{definition}
The SOC relaxation~\eqref{eqn:opf:socp} of the OPF problem \eqref{eqn:opf:nonconvex:full} is \emph{exact} if there exists a solution $(\aOpt{\mVp},\aOpt{\vp},\aOpt{\vg_{\sV}},\aOpt{\vd_{\sV}})$ in~\eqref{eqn:opf:socp} that permits a psd rank-1 completion of the partial matrix $\aOpt{\mVp}$.
\end{definition}

For such a solution, an optimizer $(\aOpt{\mV},\aOpt{\vp},\aOpt{\vg_{\sV}},\aOpt{\vd_{\sV}})$ of~\eqref{eqn:opf:sdp} with $\aOpt{\mV}={\aOpt{\vv}(\aOpt{\vv})\herm}$ is obtained that, by construction of the SDR, corresponds to an optimizer $(\aOpt{\vv},\aOpt{\vp},\aOpt{\vg_{\sV}},\aOpt{\vd_{\sV}})$ of the OPF problem~\eqref{eqn:opf:nonconvex:full}.

\begin{figure}[!t]
\centering
\includegraphics{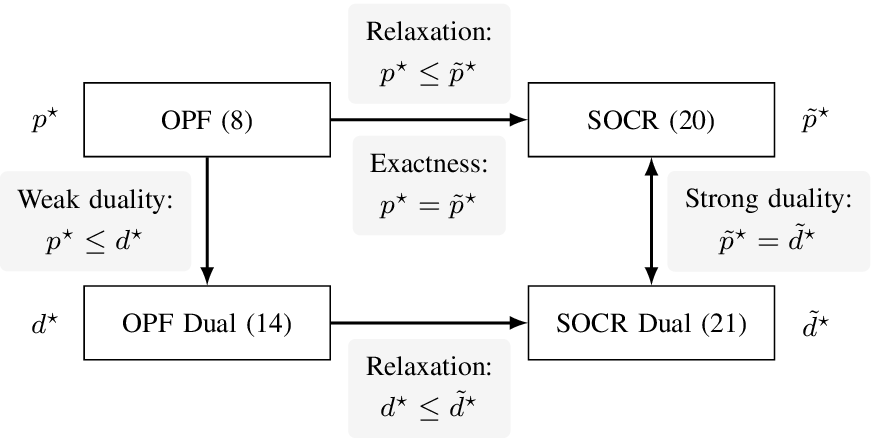}%
\caption{Interrelation of the OPF problem~\eqref{eqn:opf:nonconvex:full}, its Lagrangian dual~\eqref{eqn:opf:nonconvex:dual}, its SOC relaxation (SOCR) in~\eqref{eqn:opf:socp}, and the Lagrangian dual~\eqref{eqn:opf:socp:dual} of the SOCR. With the results in Section~\ref{sec:socpricing:method}, it follows that exactness of the relaxation implies a zero duality gap and that the optimal dual variables $\aOpt{\vLambda_{\sV}}$ in~\eqref{eqn:opf:socp:dual} are optimal in the dual~\eqref{eqn:opf:nonconvex:dual} of the OPF problem, i.e., $\aOpt{\vLambda_{\sV}}$ match the LMPs.}%
\label{fig:nodalpricing:relationsocr}%
\end{figure}

\subsection{Nodal Pricing}
	\label{sec:socpricing:method}

Analogous to SDR, the convexity of~\eqref{eqn:opf:socp} ensures that \emph{strong duality} holds, where Theorem~\ref{thm:strictfeas:socp} establishes compliance with Slater's constraint qualification.
\begin{theorem}\label{thm:strictfeas:socp}
Consider the SOC relaxation~\eqref{eqn:opf:socp} of the OPF problem~\eqref{eqn:opf:nonconvex:full}. If Assumption~\ref{ass:strictfeas} holds, there exists a feasible tuple $(\mVp,\vp,\vg_{\sV},\vd_{\sV})$ in~\eqref{eqn:opf:socp} for which $\fXi_k(\mVp)\pd\mZero$, for all~${k\in\sE}$.
\end{theorem}
\begin{IEEEproof}
This follows immediately from Theorem~\ref{thm:strictfeas:sdp} and the positive definiteness of the principal submatrices of $\mV$.
\end{IEEEproof}

However, in contrast to SDR, the Lagrangian dual of the SOC relaxation~\eqref{eqn:opf:socp} does \emph{not} match the dual~\eqref{eqn:opf:nonconvex:dual} of the OPF problem and, thus, the conclusions on LMPs in Section~\ref{sec:sdrpricing:method} do \emph{not} hold for the SOC relaxation. To show this, the \emph{dual problem} of the SOC relaxation~\eqref{eqn:opf:socp} is derived as
\begin{subequations}\label{eqn:opf:socp:dual}
\begin{align}
\aOpt{\tilde{d}} =\;&\minimize_{
	\mathclap{\hspace{2mm}\substack{\\
		\vLambda_n\in\nR^2,\,\vMu\in\nRnn^{\cM}\\[0.075em]
		\mLambda_k\in\nS^2}}}
	\quad\ 
	\vMu\tran\vb
		+ \sum_{n\in\sV} \fSigmaB_n(\vLambda_n)
		+ \sum_{n\in\sV} \fSigmaC_n(\vLambda_n)\\[0.25em]
&
\subjectto\notag\\
&\qquad
\fPsi(\vLambda_\sV,\vMu)=\sum_{k\in\sE} \mS_k\mLambda_k\mS_k\tran
	\label{eqn:opf:socp:dual:crt:psimtx}\\
&\qquad
\fpsi(\vLambda_\sV,\vMu)=\mZero\\
&\qquad
\mLambda_k\psd\mZero,
	&& \hspace{-4em} k\in\sE
	\label{eqn:opf:socp:dual:crt:lambdapsd}
\end{align}
\end{subequations}
in which~\eqref{eqn:opf:socp:crt:powbal} is dualized using $\vLambda_n\in\nR^2$, with $n\in\sV$, \eqref{eqn:opf:socp:crt:ineq} using $\vMu\in\nRnn^M$, and~\eqref{eqn:opf:socp:crt:psd} using $\mLambda_k\in\nS^2$, where $\mLambda_k\psd\mZero$ and $k\in\sE$. It can be observed that this optimization problem corresponds to the Lagrangian dual~\eqref{eqn:opf:nonconvex:dual} of the OPF problem with the psd constraint~\eqref{eqn:opf:nonconvex:dual:crt:psimtx} on $\fPsi(\vLambda_\sV,\vMu)$ being replaced by~\eqref{eqn:opf:socp:dual:crt:psimtx}. In fact, \eqref{eqn:opf:nonconvex:dual} may be regarded as a \emph{relaxation} of~\eqref{eqn:opf:socp:dual}.
\begin{lemma}\label{lem:dualfeasiblesets}
Let $(\vLambda_{\sV},\vMu,\mLambda_{\sE})$ be feasible in~\eqref{eqn:opf:socp:dual}. Then, $(\vLambda_{\sV},\vMu)$ is feasible in~\eqref{eqn:opf:nonconvex:dual}.
\end{lemma}
\begin{IEEEproof}
Due to~\eqref{eqn:opf:socp:dual:crt:lambdapsd}, it follows that $\mS_k\mLambda_k\mS_k\tran\psd\mZero$. As the sum of psd matrices is psd, \eqref{eqn:opf:socp:dual:crt:psimtx} implies $\fPsi(\vLambda_\sV,\vMu)\psd\mZero$.
\end{IEEEproof}
\begin{corollary}\label{cor:dualoptobjval}
For~\eqref{eqn:opf:nonconvex:dual} and~\eqref{eqn:opf:socp:dual}, it holds that $\aOpt{d}\leq\aOpt{\tilde{d}}$.
\end{corollary}
This observation facilitates a relation of the \emph{duality gap} of the OPF problem and \emph{exactness} of the SOC relaxation.
\begin{theorem}\label{thm:exactnesszerogap}
Assume the SOC relaxation~\eqref{eqn:opf:socp} is exact. Then, the duality gap of the OPF problem~\eqref{eqn:opf:nonconvex:full} is zero, i.e., $\aOpt{p}=\aOpt{d}$.
\end{theorem}
\begin{IEEEproof}
Exactness of the relaxation implies $\aOpt{p}=\aOpt{\tilde{p}}$ and strong duality in~\eqref{eqn:opf:socp} implies $\aOpt{\tilde{p}}=\aOpt{\tilde{d}}$. With Corollary~\ref{cor:dualoptobjval}, it follows that $\aOpt{p}\geq\aOpt{d}$. Weak duality in~\eqref{eqn:opf:nonconvex:full} implies $\aOpt{p}\leq\aOpt{d}$, thus $\aOpt{p}=\aOpt{d}$.
\end{IEEEproof}

Therefore, exactness of the SOC relaxation~\eqref{eqn:opf:socp} implies a zero duality gap of the OPF problem and, thus, qualifies the optimal dual variables $\aOpt{\vLambda_{\sV}}$ in~\eqref{eqn:opf:nonconvex:dual} as LMPs. However, this still requires the solution of the semidefinite program~\eqref{eqn:opf:nonconvex:dual}. This drawback is eliminated by the following theorem, which shows that the optimal dual variables of the SOC relaxation actually match the LMPs if the relaxation is exact.
\begin{theorem}
Let $(\aOpt{\mVp},\aOpt{\vp},\aOpt{\vg_{\sV}},\aOpt{\vd_{\sV}},\aOpt{\vLambda_{\sV}},\aOpt{\vMu},\aOpt{\mLambda_{\sE}})$ be a primal and dual optimal solution of the SOC relaxation, i.e., of~\eqref{eqn:opf:socp} and~\eqref{eqn:opf:socp:dual}, and let $\aOpt{\mVp}$ permit a psd rank-1 completion. 
Then, $(\aOpt{\vv},\aOpt{\vp},\aOpt{\vg_{\sV}},\aOpt{\vd_{\sV}},\aOpt{\vLambda_{\sV}},\aOpt{\vMu})$, where ${\aOpt{\vv}(\aOpt{\vv})\herm}$ is a psd rank-1 completion of $\aOpt{\mVp}$, is a primal and dual optimal solution of the OPF problem, i.e., of~\eqref{eqn:opf:nonconvex:full} and~\eqref{eqn:opf:nonconvex:dual}.
\end{theorem}
\begin{IEEEproof}
$(\aOpt{\vv},\aOpt{\vp},\aOpt{\vg_{\sV}},\aOpt{\vd_{\sV}})$ is optimal in~\eqref{eqn:opf:nonconvex:full} by construction of the SOC relaxation, cf. Section~\ref{sec:socpricing:relaxation}. Regarding the dual optimality, it follows from Lemma~\ref{lem:dualfeasiblesets} that $(\aOpt{\vLambda_{\sV}},\aOpt{\vMu})$ is feasible in~\eqref{eqn:opf:nonconvex:dual}. Furthermore, it follows from the proof of Theorem~\ref{thm:exactnesszerogap} that $\aOpt{d}=\aOpt{\tilde{d}}$ and, because of the equivalence of the objective functions, this implies that $(\aOpt{\vLambda_{\sV}},\aOpt{\vMu})$ is optimal in~\eqref{eqn:opf:nonconvex:dual}.
\end{IEEEproof}

Consequently, if the SOC relaxation is exact, its optimal dual variables $\aOpt{\vLambda_{\sV}}$ constitute the LMPs and its optimal primal variables $\aOpt{\vg_{\sV}}$ and $\aOpt{\vd_{\sV}}$ serve to select the bids and offers to clear the market, cf. Fig.~\ref{fig:nodalpricing:relationsocr}. Thus, under exactness, certifying the dual variables as LMPs, implementing the nodal pricing, and clearing the market is performed \emph{collectively} by solving the SOC relaxation~\eqref{eqn:opf:socp} and its dual in~\eqref{eqn:opf:socp:dual}, which is accomplished jointly and efficiently with interior-point methods.

\subsection{Bus Voltage Recovery and Relaxation Error}
	\label{sec:socpricing:recovery}

The bus voltages $\aOpt{\vv}$ associated with an optimizer $\aOpt{\mVp}$ in~\eqref{eqn:opf:socp} may be recovered with the method in~\cite[Sec.~III-B-3]{Bose2012b}. That is, the voltage magnitudes are set to the (positive) square root of the diagonal elements and the voltage angles are given by an accumulation of angle differences, which is obtained by traversing the AC subgraph and adding up the phase of the corresponding off-diagonal elements. If the SOC relaxation is exact and $\aOpt{\mVp}$ permits a {rank-1} completion, then this $\aOpt{\vv}$ is indeed optimal in the OPF problem~\eqref{eqn:opf:nonconvex:full} and ${\aOpt{\vv}(\aOpt{\vv})\herm}$ agrees with $\aOpt{\mVp}$. If the relaxation is inexact, the off-diagonal elements differ. To quantify inexactness as well as the numerical error due to finite precision solvers, it is reasonable to define the \emph{relaxation error measures}
\begin{align}
	\hspace{-3mm}
	\kappa(\aOpt{\mVp}) &= \frac{1}{\cNE}\sum_{k\in\sE}\kappa_k(\aOpt{\mVp})\\[-0.35em]
\shortintertext{and}
	\bar{\kappa}(\aOpt{\mVp}) &= \max_{k\in\sE}\;\kappa_k(\aOpt{\mVp})
\end{align}
in which the relative error related to AC branch $k\in\sE$ is
\begin{align}
	\kappa_k(\aOpt{\mVp}) =
		\left\lvert
			\frac{[\aOpt{\vv}]_{i}[\aOpt{\vv}]_{j}\conj - [\mVp]_{i,j}}{[\aOpt{\vv}]_{i}[\aOpt{\vv}]_{j}\conj}
		\right\rvert
\end{align}
where $i=\ehat(k)$ and $j=\echk(k)$.

\section{Exactness of the SOC Relaxation}
	\label{sec:exactness}

With respect to SDR-based nodal pricing in Section~\ref{sec:sdrpricing:method}, SOC-based nodal pricing can address issue (a), the practical tractability for large-scale grids. In contrast, issue (b), i.e., the characterization of exactness to promote applicability, was not yet considered and the additional relaxation even emphasizes its importance. From the discussion in Section~\ref{sec:sdrpricing:method}, it is evident that such a characterization of exactness is difficult to obtain without any structure in the grid topology. While the previous results are valid for \emph{arbitrary} hybrid transmission grids, the following study of exactness focuses on the \emph{hybrid architecture} as proposed in~\cite{Hotz2016a}, i.e., the AC subgrid is a tree while the DC subgrids may exhibit an arbitrary topology. Using these structural features, this section develops a series of results that lead to the notion of \emph{pathological price profiles}. They enable an intuitive assessment of exactness and, thus, the applicability of the SOC-based nodal pricing method. It is worth noting that \emph{no} additional technical conditions are assumed, which further distinguishes these results from the existing literature as reviewed in Section~\ref{sec:sdrpricing:method}.

\subsection{Sufficient Condition for Exactness}

By virtue of the \emph{tree topology} of the AC subgrid, it follows from~\cite[Th.~5]{Bose2012b} that there exists a (unique) psd rank-1 completion of an optimizer $\aOpt{\mVp}$ in~\eqref{eqn:opf:socp} if
\begin{equation}\label{eqn:suffcond:socp}
	\rank\big(\fXi_k(\aOpt{\mVp})\big) = 1,
	\qquad
	\forall k\in\sE.
\end{equation}
On the other hand, strong duality implies that every primal optimal $\aOpt{\mVp}$ in~\eqref{eqn:opf:socp} and dual optimal $\aOpt{\mLambda_\sE}$ in~\eqref{eqn:opf:socp:dual} satisfy the \emph{complementary slackness}
\begin{align}\label{eqn:opf:socp:compslack}
	\trace\big(\aOpt{\mLambda_k}\,\fXi_k(\aOpt{\mVp})\big) = 0,
	\qquad
	\forall k\in\sE\,.
\end{align}
By the following theorem, the necessary condition for \emph{optimality} in~\eqref{eqn:opf:socp:compslack} is related to the sufficient condition for \emph{exactness} in~\eqref{eqn:suffcond:socp}. To this end, the sets $\sA_k$, with $k\in\sE$, are defined~as
\begin{align}\label{eqn:akdef}
	\hspace{-0.4em}
	\sA_k = \big\{ (\vLambda_\sV,\vMu) \in \nR^2_\sV\times\nR^\cM\!:
		[\fPsi(\vLambda_\sV,\vMu)]_{\ehat(k),\echk(k)} = 0 \hskip0.75pt\big\}.
\end{align}
\begin{theorem}\label{thm:exactness:notinak}
Let $(\aOpt{\mVp},\aOpt{\vp},\aOpt{\vg_\sV},\aOpt{\vd_\sV},\aOpt{\vLambda_\sV},\aOpt{\vMu},\aOpt{\mLambda_\sE})$ be an optimizer of the primal and dual problem~\eqref{eqn:opf:socp} and~\eqref{eqn:opf:socp:dual}. For any $k\in\sE$, if $(\aOpt{\vLambda_\sV},\aOpt{\vMu}) \notin \sA_k$, then $\rank\big(\fXi_k(\aOpt{\mVp})\big)=1$.
\end{theorem}
\begin{IEEEproof}
See Appendix~\ref{apx:thm:exactness:notinak}.
\end{IEEEproof}
\begin{corollary}\label{cor:exactness:suffcond}
If $(\aOpt{\vLambda_\sV},\aOpt{\vMu}) \notin \bigcup_{k\in\sE} \sA_k$, then, for all $k\in\sE$, $\rank\big(\fXi_k(\aOpt{\mVp})\big)=1$ and the SOC relaxation is \emph{exact}.
\end{corollary}

The sufficient condition for \emph{exactness} in~\eqref{eqn:suffcond:socp} thus translates to the \emph{avoidance} of the sets $\sA_k$ in dual optimality, which are linear subspaces as established by the following theorem.
\begin{theorem}\label{thm:aksubspace}
For all $k\in\sE$, $\sA_k$ is a proper linear subspace and $\dim(\nR^2_\sV\times\nR^\cM)-2 \leq \dim(\sA_k) \leq \dim(\nR^2_\sV\times\nR^\cM)-1$.
\end{theorem}
\begin{IEEEproof}
See Appendix~\ref{apx:thm:aksubspace}.
\end{IEEEproof}
Indeed, if the dual optimal variables $(\aOpt{\vLambda_\sV},\aOpt{\vMu})$ lie in the subspace $\sA_k$, the dual variable $\aOpt{\mLambda_k}$ of the corresponding constraint in~\eqref{eqn:opf:socp:crt:psd} is zero as shown below. This indicates potential inactivity of the constraint, in which case the associated {$2$$\times$$2$} principal submatrix of $\aOpt{\mVp}$ is rank-2.
\begin{theorem}\label{thm:lambdakzero}
Let $(\aOpt{\vLambda_\sV},\aOpt{\vMu},\aOpt{\mLambda_\sE})$ be an optimizer of~\eqref{eqn:opf:socp:dual}. If $(\aOpt{\vLambda_\sV},\aOpt{\vMu})\in\sA_k$, then $\aOpt{\mLambda_k}=\mZero$.
\end{theorem}
\begin{IEEEproof}
See Appendix~\ref{apx:thm:lambdakzero}.
\end{IEEEproof}

With Theorem~\ref{thm:aksubspace} and~\ref{thm:lambdakzero}, the \emph{sufficient condition} for \emph{exactness} in Corollary~\ref{cor:exactness:suffcond} states that exactness of the relaxation may only be lost if the optimal dual variables $(\aOpt{\vLambda_\sV},\aOpt{\vMu})$ lie in a union of $\cNE$ subspaces in an $2\cNV+\cM$ dimensional space, where $\cNE=\cNV-1$ as the AC subgrid is a tree (cf.~\cite[Corollary~1]{Hotz2016a}).

\subsection{Relation of Exactness and Locational Marginal Prices}

The results in the previous section provide an insight into exactness on a rather abstract level. To render them intuitively accessible, they shall be related to LMPs.
\begin{theorem}\label{thm:exactness:pospriceacb}
Let $(\aOpt{\vLambda_\sV},\aOpt{\vMu},\aOpt{\mLambda_\sE})$ be an optimizer of~\eqref{eqn:opf:socp:dual}. Consider any $k\in\sE$ and let $i=\ehat(k)$ and $j=\echk(k)$. If $\aOpt{\vLambda_i}\geq\mZero$, $\aOpt{\vLambda_j}\geq\mZero$, and $[\aOpt{\vLambda_i}]_1+[\aOpt{\vLambda_j}]_1>0$, then $(\aOpt{\vLambda_\sV},\aOpt{\vMu}) \notin \sA_k$.
\end{theorem}
\begin{IEEEproof}
See Appendix~\ref{apx:thm:exactness:pospriceacb}.
\end{IEEEproof}
\begin{corollary}
If $[\aOpt{\vLambda_n}]_1>0$ and $[\aOpt{\vLambda_n}]_2\geq 0$, for all $n\in\sV$,~then $\rank\big(\fXi_k(\aOpt{\mVp})\big)=1$, for all $k\in\sE$, and the relaxation is \emph{exact}.
\end{corollary}

Considering that $\aOpt{\vLambda_n}$ is an LMP-based \emph{price vector}, this result gives rise to an interesting observation: If the LMP of active power is positive and the LMP of reactive power is nonnegative at all buses, exactness of the relaxation is guaranteed. Indeed, typically at the majority of buses such a price structure is observed, cf. the discussion in Appendix~\ref{apx:lmpimplications}. While this appears evident for active power due to convex and increasing generation cost functions, it may not be as obvious for reactive power, which is usually associated with zero cost and benefit. In this case, if the constraints on reactive power are not binding, its LMP is zero and, if they are binding, it is mostly due to a demand for capacitive reactive power that leads to a positive price. However, buses with negative LMPs do arise, e.g., in case of demand for more inductive reactive power. Still, exactness obtains as long as the optimal dual variables are not forced into the union of the subspaces $\sA_k$, cf. Theorem~\ref{thm:exactness:notinak} and Corollary~\ref{cor:exactness:suffcond}. This unfortunate constellation that may compromise exactness is termed as follows.
\begin{definition}
If an optimizer $(\aOpt{\vLambda_\sV},\aOpt{\vMu},\aOpt{\mLambda_\sE})$ of~\eqref{eqn:opf:socp:dual} satisfies
\begin{equation}\label{eqn:suffcond:pathological}
(\aOpt{\vLambda_\sV},\aOpt{\vMu})\in\bigcup_{k\in\sE} \sA_k
\end{equation}
then the $\cNV$-tuple $\aOpt{\vLambda_\sV}$ is a \emph{pathological price profile}.
\end{definition}

\begin{figure}[!t]
\centering
\includegraphics{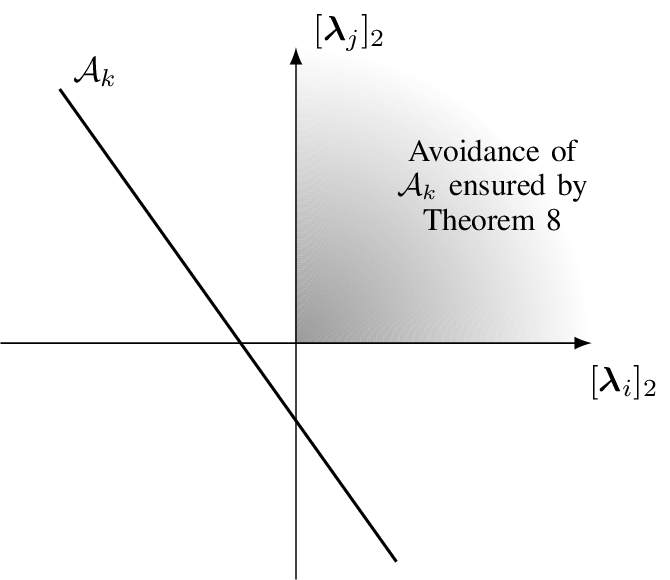}%
\caption{Qualitative illustration of some pathological price profiles: Exactness of the SOC relaxation may only be lost if the optimal dual variables $[\aOpt{\vLambda_i}]_2$ and $[\aOpt{\vLambda_j}]_2$, which constitute a point in this plane, happen to lie on the \mbox{$\sA_k$-line}. In all other cases, exactness is guaranteed and $[\aOpt{\vLambda_i}]_2$ and $[\aOpt{\vLambda_j}]_2$ equal the LMP for reactive power at bus $i$ and $j$, respectively.}%
\label{fig:exactness:subspace}%
\end{figure}

The particular value of~\eqref{eqn:suffcond:pathological} is its \emph{characterization of exactness}. To illustrate this, let $(\aOpt{\vLambda_\sV},\aOpt{\vMu},\aOpt{\mLambda_\sE})$ be an optimizer of~\eqref{eqn:opf:socp:dual} and assume that all subspaces except $\sA_k$ are avoided, i.e., $(\aOpt{\vLambda_\sV},\aOpt{\vMu})\not\in\bigcup_{r\in\sE\setminus\{k\}} \sA_r$. In this case, exactness is ensured if $\sA_k$ is avoided, i.e., if $(\aOpt{\vLambda_\sV},\aOpt{\vMu})\not\in\sA_k$, while it is potentially lost otherwise. Let $i=\ehat(k)$ and $j=\echk(k)$ denote the adjacent buses of AC branch $k$. Assume the LMP for \emph{active} power at bus $i$ and $j$ is positive and $[\aOpt{\vLambda_i}]_1+[\aOpt{\vLambda_j}]_1>0$. In this case, exactness of the relaxation is determined by the LMP for \emph{reactive} power at bus $i$ and $j$, i.e., $[\aOpt{\vLambda_i}]_2$ and $[\aOpt{\vLambda_j}]_2$. The potential loss of exactness can be visualized qualitatively as shown in Fig.~\ref{fig:exactness:subspace}.\footnote{Figure~\ref{fig:exactness:subspace} is the projection of a cut through the ambient space of $\sA_k$. This cut is parallel to the $[\vLambda_i]_2$-$[\vLambda_j]_2$-plane and includes $(\aOpt{\vLambda_\sV},\aOpt{\vMu})$. The shape of $\sA_k$ follows from Theorem~\ref{thm:aksubspace}. In this respect, the cut illustrates the worst case, i.e., a line, see also~\eqref{eq:ak:nullspace}. The location of $\sA_k$ follows from Theorem~\ref{thm:exactness:pospriceacb}, which states that it does not intersect the nonnegative quadrant.} Exactness may only be lost if $[\aOpt{\vLambda_i}]_2$ and $[\aOpt{\vLambda_j}]_2$ combine to a point on the $\sA_k$-line. Intuitively, this appears unlikely, or \emph{pathological}, considering the nature of LMPs as well as their sensitivity to operating conditions. The argument of this example extends similarly to further subspaces and LMPs.\footnote{For every AC branch $k$, avoidance of the associated subspace $\sA_k$ depends \emph{exclusively} on the LMPs at the adjacent buses. This follows from $\sA_k$ in~\eqref{eqn:akdef} and $\fPsi(\vLambda_\sV,\vMu)$ in~\eqref{eqn:psimtxdef} as well as the structure of the constraint matrices documented in~\cite[Appendix~B]{Hotz2016a}. The element in row $\ehat(k)$ and column $\echk(k)$ may only be nonzero for constraint matrices related to AC branch $k$, i.e., only the dual variables of these constraints determine the avoidance of $\sA_k$.} These insights motivate the conjecture that the price profiles in (28) are pathological and, as a consequence, that the relaxation is typically exact for the hybrid architecture. This is supported by the following simulation results.

\begin{figure*}[!t]
\centering
\includegraphics{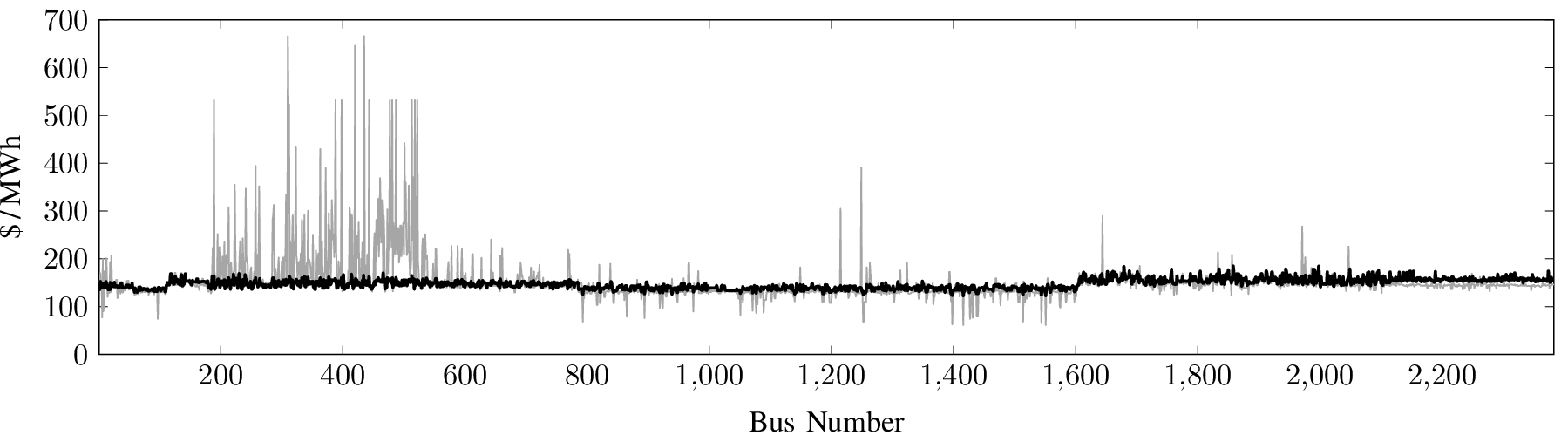}%
\vskip-0.5em%
\caption{Case Study 1: LMP profile of active power (with respect to the system buses) for the HTG (black line) and the ACG (gray line).}%
\label{fig:sim:cs1:lmpprofile}%
\end{figure*}

\begin{figure}[!t]
\centering
\includegraphics{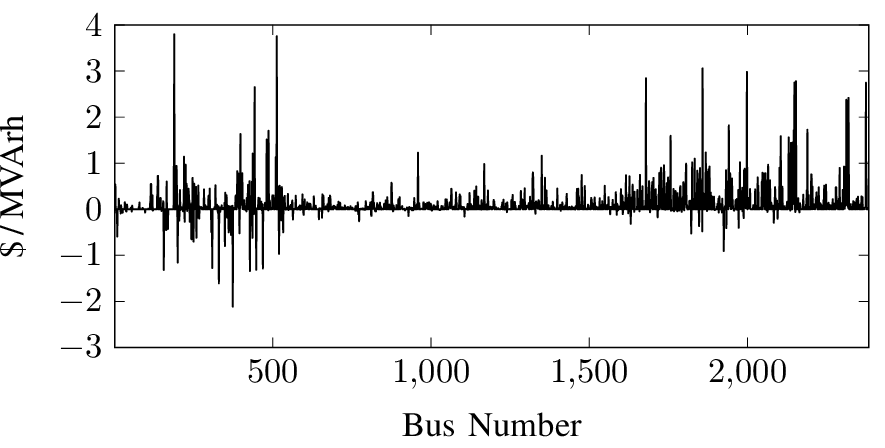}%
\vskip-0.5em%
\caption{Case Study 1: LMP profile of reactive power for the HTG.}
\label{fig:sim:cs1:lmpreactive}
\end{figure}

\section{Simulation Results}
	\label{sec:simulation}

In the following, the proposed nodal pricing method is utilized to highlight benefits of the hybrid architecture. To this end, the 2383-bus test case ``{case2383wp.m}'' provided by the power system simulation package \textsc{Matpower}~\cite{Zimmerman2016a,Zimmerman2011a} is considered, which represents the Polish transmission grid during winter peak conditions. It is converted to a hybrid transmission grid via a proposed upgrade strategy. Subsequently, the original AC transmission grid (ACG) and the hybrid transmission grid (HTG) are compared in four case studies. For the HTG, the LMPs are determined using the nodal pricing method presented in Section~\ref{sec:socpricing:method}, where the SOC relaxation~\eqref{eqn:opf:socp} and its dual~\eqref{eqn:opf:socp:dual} are solved simultaneously with \textsc{Mosek}~\cite{MOSEK-ApS2016a},\footnote{Default settings with \texttt{MSK\_DPAR\_INTPNT\_CO\_TOL\_PFEAS\;=\;}$5\cdot 10^{-8}$.} an interior-point optimizer for large-scale conic optimization problems. For the ACG, approximate LMPs are determined via a DC\;OPF, as the SOC relaxation proves to be inexact. Thus, the presented LMPs of the ACG are given by the dual variables of the power balance constraints in the ``DC power flow'' model, which is a linearization of the AC power flow that ignores losses as well as reactive power and presumes a perfectly flat voltage profile. The DC\;OPF for the ACG is solved with \textsc{Matpower}~\cite{Zimmerman2016a} using \textsc{Mosek}~\cite{MOSEK-ApS2016a} as the solver back end.

\subsection{Preprocessing}

The test case is preprocessed to support the case studies.

\subsubsection{System-Restricting Branches}

If the load is increased by more than $5.0\,$\% the DC\;OPF of the ACG becomes infeasible. The system is actually limited by only two out of $2896$ AC branches, i.e., branch $2239$ and $2862$. To enable case studies under high load, their line rating is increased by $35\,$\%.

\subsubsection{Parallel Branches}

The test case contains $10$ pairs of parallel AC branches. The parallel branches $18$ and $19$ exhibit the same line parameters but different transformer tap ratios, i.e., $1.0646$ and $1.0523$. For maximum capacity, the tap ratio of branch $19$ is set to $1.0646$. Parallel AC branches are then combined into an equivalent single AC branch to comply with the hybrid transmission grid model (cf.~\cite[Def.~5]{Hotz2016a}).

\subsubsection{Lossless Branches}

The test case contains $195$ lossless AC branches. The hybrid transmission grid model assumes lossy AC branches, cf.~\cite[Def.~7 and~8]{Hotz2016a}. To this end, all lossless branches are imposed with negligible losses by setting their series resistance to $10^{-5}\,$p.u..

\subsection{Upgrade to the Hybrid Architecture}
	\label{sec:simulation:upgrade}

Possible options for the upgrade to the hybrid architecture correspond to spanning trees of the ACG, see~\cite{Hotz2016a}. Their number is given by Kirchhoff's matrix tree theorem, cf. e.g.~\cite{Merris1994a}. In the numerical computing environment \textsc{Matlab}~\cite{MATLAB2016a}, its evaluation results in an overflow, implying that there are more than $10^{308}$ spanning trees. Thus, the evaluation of all upgrade options as in~\cite{Hotz2016a} is not tractable and, on that account, we propose a heuristic upgrade strategy. To this end, every AC branch is associated with an \emph{upgrade suitability measure} to obtain a weighted graph. Therein, the minimum spanning tree is identified and all AC branches outside this tree are upgraded to DC branches. As long-distance transmission is a common application of HVDC lines~\cite{Bahrman2007a}, long transmission lines are favored for the upgrade by using the series resistance of an AC branch as its upgrade suitability measure. Therewith, $17.46\,$\% of all AC branches are selected and converted to (bidirectional) HVDC lines with an exemplary loss factor of $3.5\%$ and a capacity that coincides with the AC line rating to focus on the influence of the architecture. The reactive power capability is selected conservatively as $25\,$\% of the active power capacity, cf. e.g.~\cite{ABB-Grid-Systems2012a}.

\subsection{Case Study 1: LMP Profile Equalization}

In this test case, the load totals to $24.6\,$GW, i.e., $82.99\,$\% of the active power generation capacity (P-capacity). The corresponding LMPs for active power are depicted in Fig.~\ref{fig:sim:cs1:lmpprofile}. For the ACG, the LMP fluctuates excessively between $61.40\,$\$/MWh and $665.69\,$\$/MWh. Furthermore, a subsequent AC power flow reveals that the model mismatch necessitates substantial $667.3\,$MW of slack power at the reference bus. In contrast, the LMP profile for the HTG is remarkably flat and all LMPs lie in the narrow interval from $122.82\,$\$/MWh to $184.71\,$\$/MWh, while \emph{no} slack power is required. Considering that the generation cost functions are convex and monotonically increasing, it follows from the definition of LMPs that the disappearing of downward peaks in the LMP profile implies an improved utilization of generators with low marginal costs, while the vanishing of upward peaks states that loads obtain access to less expensive generation. This ``equalization'' of the LMP profile thus implies that restrictions on trades are reduced and, hence, the separation of nodal power markets due to grid limitations is mitigated.

\begin{table}[!t]
\renewcommand{\arraystretch}{\tabstretch}
\caption{Relaxation Error for the HTG}
\label{tab:sim:recerr}
\centering
\includegraphics{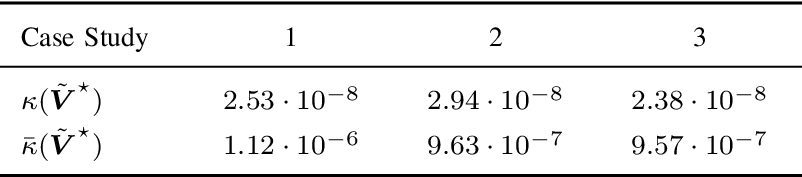}%
\end{table}

For the HTG, the accuracy of LMPs is established by the proposed nodal pricing method, while the equalization arises from the flexibility induced by the hybrid architecture. This is verified by performing a DC\;OPF for the HTG, where the equalization is also observed in these approximate LMPs, which range from $127.36\,$\$/MWh to $158.50\,$\$/MWh.

Table~\ref{tab:sim:recerr} shows that the SOC relaxation is exact for the HTG. In fact, the LMP of active power is always positive and Fig.~\ref{fig:sim:cs1:lmpreactive} illustrates that the LMP of reactive power is nonnegative at the majority of buses. Accordingly, Theorem~\ref{thm:exactness:pospriceacb} ensures the avoidance of $52.23\,$\% of the critical subspaces $\sA_k$. Furthermore, the smallest magnitude of the off-diagonal elements of $\fPsi(\aOpt{\vLambda_\sV},\aOpt{\vMu})$ on the AC subgraph is $9.659$, i.e., the remaining critical subspaces in~\eqref{eqn:akdef} are also avoided and Corollary~\ref{cor:exactness:suffcond} guarantees exactness. This result supports the conjecture that it is unlikely to match the critical subspaces in~\eqref{eqn:suffcond:pathological} and, thus, to observe a pathological price profile, cf. Fig.~\ref{fig:exactness:subspace}. Case Study~$4$ substantiates this conclusion later on. Finally, it shall be noted that the computational effort is rather moderate with the solver time of \textsc{Mosek} being less than $2$~seconds on a standard office notebook.

If the proposed nodal pricing method is applied to the ACG, it is observed that the relaxation is inexact and, thus, the optimal dual variables $\aOpt{\vLambda_{\sV}}$ do \emph{not} constitute viable nodal prices. The relaxation error amounts to $\kappa(\aOpt{\mVp})=9.80\cdot 10^{-3}$ and $\bar{\kappa}(\aOpt{\mVp})=0.472$. After bus voltage recovery, the average error in the apparent power balance is $1221\,$MVA, i.e., the system model is disregarded and the results are inapplicable.

\subsection{Case Study 2: Nodal Pricing at High Load}

\begin{figure}[!t]
\centering
\includegraphics{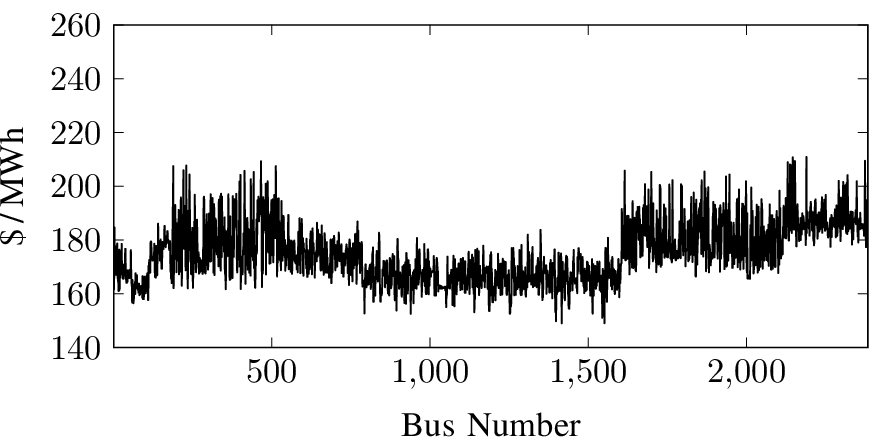}%
\vskip-0.5em%
\caption{Case Study 2: LMP profile of active power for the HTG at high load.}
\label{fig:sim:cs2:lmpprofile}
\end{figure}

This case study investigates the impact of high load on the LMP profile. In the ACG, the load can be scaled up by $9.9\,$\% until \textsc{Matpower} reports infeasibility. At this load, which amounts to $91.20\,$\% of P-capacity, the LMPs of the ACG fluctuate vastly between $-3{,}229.60\,$\$/MWh and $15{,}506.91\,$\$/MWh. In contrast, the HTG supports a load increase of $16.8\,$\%, totaling the load to $96.92\,$\% of P-capacity. Considering that the transmission losses are $3.08\,$\% of P-capacity, it follows that \emph{all} generation can be utilized. Figure~\ref{fig:sim:cs2:lmpprofile} illustrates that even under these extreme conditions the LMP profile remains reasonably flat. Thus, the hybrid architecture can not only increase the effective capacity, but also maintain adequate nodal prices under very high load.

\subsection{Case Study 3: Loadability After Generation Expansion}

\begin{figure}[!t]
\centering
\includegraphics{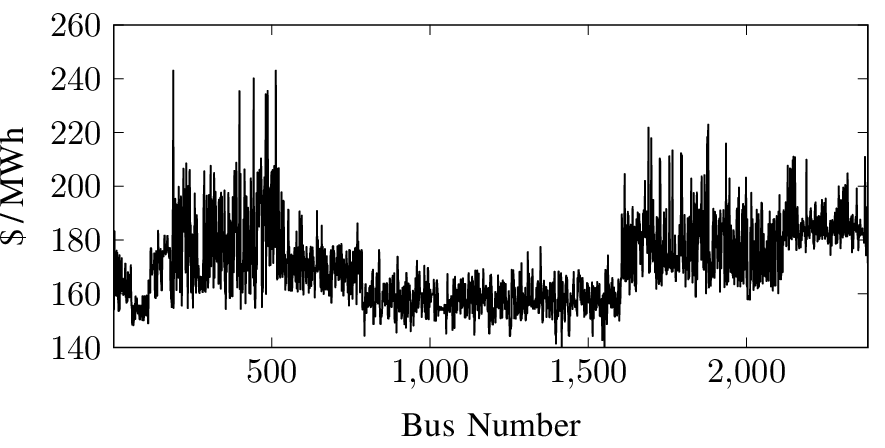}%
\vskip-0.5em%
\caption{Case Study 3: LMP profile of active power for the HTG at very high load after expanding generation by $15\,\%$.}
\label{fig:sim:cs3:lmpprofile}
\end{figure}

To further investigate loadability, the generation capacity is increased by $15\,\%$. Under these conditions, the ACG supports a load increase of $16.7\,$\%, which totals the load to $84.21\,$\% of P-capacity and results in LMPs fluctuating between $-28{,}241.51\,$\$/MWh and $69{,}620.31\,$\$/MWh. In contrast, the HTG supports a substantially higher load increase of $29.7\,$\% that totals the load to $93.59\,$\% of P-capacity. With transmission losses of $3.31\,$\% of P-capacity, it follows that still almost all generation can be utilized and Fig.~\ref{fig:sim:cs3:lmpprofile} shows that the HTG maintains a reasonable LMP profile. The effective capacity of the HTG is thus $11.14\,$\% higher compared to the ACG and, additionally, it offers the major advantage of desirable trading conditions. Concluding, it should be emphasized that this capacity gain arises \emph{exclusively} from the flexibility induced by the hybrid architecture, as the line capacity is not uprated during the conversion of AC lines to HVDC in Section~\ref{sec:simulation:upgrade}.

\subsection{Case Study 4: Exactness of the SOC Relaxation}

\begin{figure}[!t]
\centering
\includegraphics{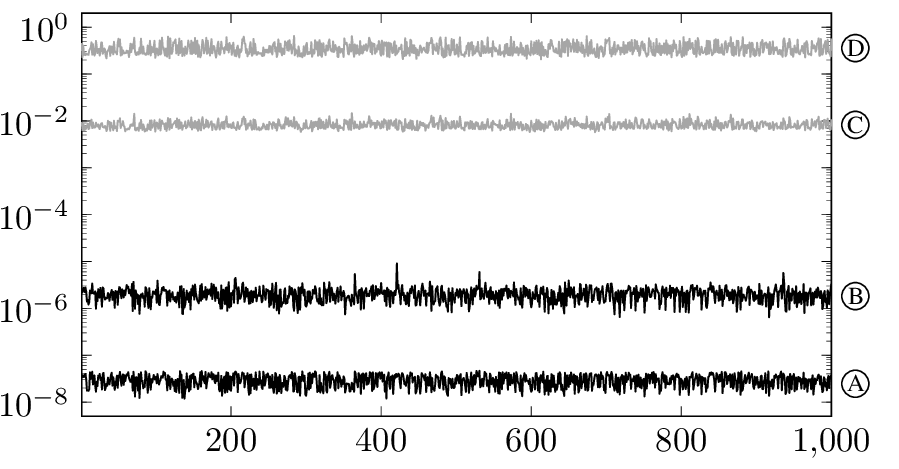}%
\vskip-0.25em%
\caption{Case Study 4: Relaxation error for $1000$ scenarios with different operating conditions. For the HTG, $\kappa(\aOpt{\mVp})$ and $\bar{\kappa}(\aOpt{\mVp})$ is given by $\protect\circledA$ and~$\protect\circledB$, respectively. For the ACG, $\kappa(\aOpt{\mVp})$ and $\bar{\kappa}(\aOpt{\mVp})$ is given by $\protect\circledC$ and~$\protect\circledD$.}
\label{fig:sim:cs4:kappa}
\end{figure}

Above case studies utilize the proposed nodal pricing method to accurately identify the LMPs of the HTG. The nodal pricing method requires exactness of the relaxation, which is shown to obtain in these cases, i.e., pathological price profiles are not observed. Case Study~$1$ provides some insights on LMPs to motivate the unlikeliness of pathological price profiles under normal operating conditions. In the following, this rationale is substantiated by analyzing exactness under variations of load and generation cost. To this end, the load at \emph{every} bus is varied randomly between $25\,\%$ and $125\,\%$ by scaling with a random factor sampled from a uniform distribution on $(0.25,1.25)$. Analogously, the cost of \emph{every} generator is varied randomly between half and twice the amount by scaling the generation cost function with a random factor sampled from a uniform distribution on $(0.5,2)$. This perturbation of load and generation cost is repeated to generate $1000$ scenarios with varying operating conditions. Fig.~\ref{fig:sim:cs4:kappa} depicts the relaxation error of the HTG as well as the ACG for these scenarios. Throughout, the ACG exhibits a substantial relaxation error, rendering the nodal pricing method inapplicable. In contrast, the HTG consistently obtains exactness, which confirms that pathological price profiles are indeed unlikely.

\section{Conclusion}
	\label{sec:conclusion}

This paper identified advantages of the recently proposed hybrid transmission grid architecture in a market context, which are twofold. On one hand, this hybrid architecture was shown to provide computational benefits. To this end, an accurate and efficient nodal pricing method based on a second-order cone relaxation of the optimal power flow problem was presented, whose applicability requires exactness of the relaxation. By means of the Polish transmission grid, it was illustrated that exactness is typically not obtained for a conventional grid topology. In contrast, the hybrid architecture ensures exactness as long as the locational marginal prices do not coincide with certain pathological price profiles, which were shown to be unlikely under normal operating conditions. Therewith, the hybrid architecture enables accurate and computationally efficient nodal pricing. On the other hand, this system structure was shown to introduce operational benefits. To this end, it was utilized that locational marginal prices reveal operational characteristics of the grid. Therewith, it was illustrated that implementing capacity expansion via the transition to the hybrid architecture can substantially reduce grid-induced restrictions on power trades and improve the utilization of generation facilities.

\appendices

\begin{figure*}[!t]
\centering%
\captionsetup[subfloat]{captionskip=-0.75ex}
\subfloat[]{\hspace{-3.5mm}%
\includegraphics{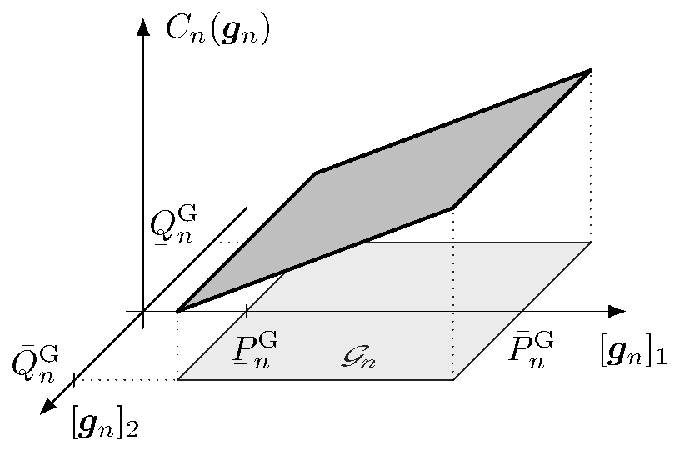}%
\label{fig:cost:3d}}%
\hspace{4.5mm}%
\subfloat[]{%
\includegraphics{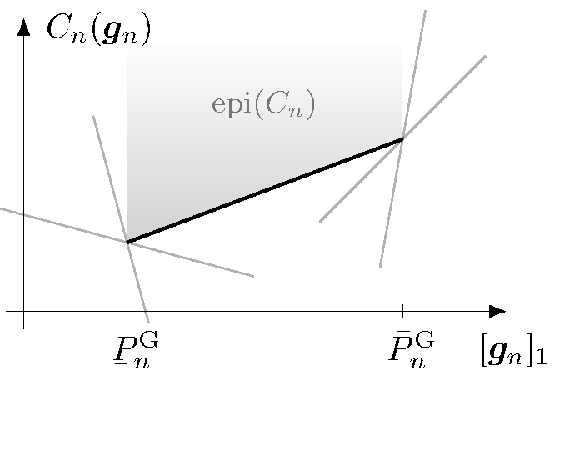}%
\label{fig:cost:active}}%
\hspace{2.5mm}%
\subfloat[]{%
\includegraphics{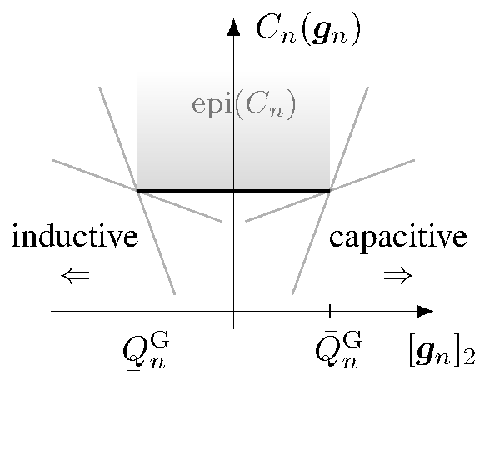}%
\label{fig:cost:reactive}\hspace{-3mm}}%
\caption{Illustration of the generation cost function $\fC_n:\sG_n\rightarrow\nR$ in~\eqref{eq:cost:cn}, where (a) shows the function over its entire domain, (b) shows a cut along the first dimension, and (b) shows a cut along the second dimension. The gray lines in (b) and (c) depict exemplary supporting hyperplanes of the epigraph. Their slope constitutes a subgradient at the respective point and in the respective dimension. The set of all subgradients is the subdifferential at that point.}
\label{fig:cost}
\end{figure*}

\section{Proof of Theorem~\ref{thm:strictfeas:sdp}}
	\label{apx:thm:strictfeas:sdp}

Assumption~\ref{ass:strictfeas} states that there exists a feasible tuple $(\vv,\vp,\vg'_{\sV},\vd'_{\sV})$ in~\eqref{eqn:opf:nonconvex:compact} for which $\mC\mkern-0.5mu\vec(\vv\vv\herm) + \aBar{\mC}\vp < \vb$ and $\vg'_n-\vd'_n\in\interior(\sG_n-\sD_n)$, for all $n\in\sV$. Let $\mV = \vv\vv\herm + \varepsilon\mI\,$, where $\varepsilon>0$ is an arbitrary positive scalar and $\mI$ is the $\cNV$$\times$$\cNV$ identity matrix, i.e., $\mV\pd\mZero$. Therewith,
\begin{align}
	&\mC\mkern-0.5mu\vec(\mV) + \aBar{\mC}\vp
		= \mC\mkern-0.5mu\vec(\vv\vv\herm) + \aBar{\mC}\vp + \varepsilon\mC\mkern-0.5mu\vec(\mI)
\shortintertext{and, with~\eqref{eqn:opf:nonconvex:compact:powbal},}
	&\mB_n\mkern-1mu\vec(\mV) + \aBar{\mB}_n\vp
		= \vg'_n - \vd'_n + \varepsilon\mB_n\mkern-1mu\vec(\mI)
\end{align}
for all $n\in\sV$. For $\varepsilon$ sufficiently small, there exist $\vg_n\in\sG_n$ and $\vd_n\in\sD_n$ such that
\begin{equation}
	\vg_n - \vd_n = \vg'_n - \vd'_n + \varepsilon\mB_n\mkern-1mu\vec(\mI)\,.
\end{equation}
for all $n\in\sV$. Consequently, there exists some $\varepsilon>0$ such that $(\mV,\vp,\vg_{\sV},\vd_{\sV})$ is feasible in~\eqref{eqn:opf:sdp}.

\section{Proof of Theorem~\ref{thm:exactness:notinak}}
	\label{apx:thm:exactness:notinak}

For all $n\in\sV$, the lower bound in~\eqref{eqn:model:crt:vm} on the bus voltage magnitude translates to a constraint
\begin{equation}\label{eqn:model:crt:vm:lb:socp}
	\trace(\mM_n\mVp)\geq\aLB{V_n}^2
	\ \ \Leftrightarrow\ \
	\vec(-\mM_n\tran)\tran\vec(\mVp)\leq-\aLB{V_n}^2
\end{equation}
in~\eqref{eqn:opf:socp} and, as $\mM_n = \ve_n\ve_n\tran$ (cf.~\cite[Sec.~III]{Hotz2016a}) and $\aLB{V_n}>0$, this ensures strictly positive diagonal elements in $\mVp$ and, thus, $\rank\big(\fXi_k(\aOpt{\mVp})\big)\geq 1$. By the absence of (anti-) parallel AC branches (cf.~\cite[Def.~5]{Hotz2016a}), it follows from~\eqref{eqn:opf:socp:dual:crt:psimtx} that $[\aOpt{\mLambda_k}]_{1,2} = [\fPsi(\aOpt{\vLambda_\sV},\aOpt{\vMu})]_{\ehat(k),\echk(k)}$, thus $(\aOpt{\vLambda_\sV},\aOpt{\vMu}) \notin \sA_k$ implies ${[\aOpt{\mLambda_k}]_{1,2}\neq 0}$ and $\rank(\aOpt{\mLambda_k})\geq 1$. In conjunction with~\eqref{eqn:opf:socp:crt:psd}, \eqref{eqn:opf:socp:dual:crt:lambdapsd}, \eqref{eqn:opf:socp:compslack}, and the rank-nullity theorem, this yields $\rank\big(\fXi_k(\aOpt{\mVp})\big)\leq 1$, cf.~\cite[Lemma~2]{Hotz2016a}. Hence, $\rank\big(\fXi_k(\aOpt{\mVp})\big)=1$.

\section{Proof of Theorem~\ref{thm:aksubspace}}
	\label{apx:thm:aksubspace}

Let $\vXi=[\vLambda_1\tran,\ldots,\vLambda_{\cNV}\tran,\vMu\tran]\tran\in\nR^{2\cNV+\cM}$ be a stacking of the tuple $(\vLambda_\sV,\vMu)$. Therewith, $\sA_k$ can be rewritten as
\begin{align}
	\sA_k
		&= \{\, \vXi\in\nR^{2\cNV+\cM} : \mL_k\vXi=\mZero \,\}
		 = \nullspace(\mL_k)
		\label{eq:ak:nullspace}
\intertext{where $\mL_k\in\nR^{2\times(2\cNV+\cM)}$ is}
	\mL_k &=
		\begin{bmatrix}
			\real(\mL'_{k,1}) & \ldots & \real(\mL'_{k,\cNV}) & \real(\mL_k'') \\
			\imag(\mL'_{k,1}) & \ldots & \imag(\mL'_{k,\cNV}) & \imag(\mL_k'') \\
		\end{bmatrix}.
		\label{eq:ak:lk}
\end{align}
Therein, $\mL'_{k,n}=\big[[\mP_n]_{i,j}\,,\,[\mQ_n]_{i,j}\big]\in\nC^{1\times 2}$ and $\mL_k''=\big[[\mC_1]_{i,j}\,,\,\ldots\,,\,[\mC_{\cM}]_{i,j}\big]\in\nC^{1\times \cM}$, with $i=\ehat(k)$ and $j=\echk(k)$. It follows from~\cite[Appendix~B, eq.~(8), and Corollary~2]{Hotz2016a} that $[\mP_i]_{i,j}\neq 0$, thus $\mL_k\neq\mZero$ and $1\leq\rank(\mL_k)\leq 2$. The rank-nullity theorem completes the proof.

\section{Proof of Theorem~\ref{thm:lambdakzero}}
	\label{apx:thm:lambdakzero}

It follows from $(\aOpt{\vLambda_\sV},\aOpt{\vMu})\in\sA_k$ that
\begin{equation}
	[\aOpt{\mLambda_k}]_{1,2}
		= [\aOpt{\mLambda_k}]_{2,1}\conj
		= [\fPsi(\aOpt{\vLambda_\sV},\aOpt{\vMu})]_{\ehat(k),\echk(k)}
		= 0
\end{equation}
cf. Appendix~\ref{apx:thm:exactness:notinak}. Moreover, \eqref{eqn:opf:socp:dual:crt:lambdapsd} implies $[\aOpt{\mLambda_k}]_{1,1},[\aOpt{\mLambda_k}]_{2,2}\geq 0$ and from \eqref{eqn:model:crt:vm:lb:socp} it follows that $[\fXi_k(\aOpt{\mVp})]_{1,1}\,,\,[\fXi_k(\aOpt{\mVp})]_{2,2}>0$. Therewith, \eqref{eqn:opf:socp:compslack} implies $[\aOpt{\mLambda_k}]_{1,1}=[\aOpt{\mLambda_k}]_{2,2}=0$, as it requires the summation of nonnegative terms to zero. Hence, $\aOpt{\mLambda_k}=\mZero$.

\section{Proof of Theorem~\ref{thm:exactness:pospriceacb}}
	\label{apx:thm:exactness:pospriceacb}

Due to the nonnegativity of $\aOpt{\vLambda_\sV}$ and $\aOpt{\vMu}$, it follows from~\cite[Th.~2 and Th.~3]{Hotz2016a} that $[\fPsi(\aOpt{\vLambda_\sV},\aOpt{\vMu})]_{i,j}\in\sH_k$, where $\sH_k\subset\nC$ is a half-space with normal $\rho_k\in\nC\setminus\{0\}$. Furthermore,~\cite[Appendix~B, eq.~(8), and Corollary~2]{Hotz2016a} implies
\begin{equation}
	\real(\rho_k\conj [\mP_i]_{i,j})
		= \real(\rho_k\conj [\mP_j]_{i,j})
		< 0\,.
\end{equation}
Therefore, \cite[Def.~16]{Hotz2016a} states $[\mP_i]_{i,j},[\mP_j]_{i,j}\in\interior(\sH_k)$. With $\aOpt{\vLambda_i},\aOpt{\vLambda_j}\geq\mZero$ and $[\aOpt{\vLambda_i}]_1+[\aOpt{\vLambda_j}]_1>0$, \cite[Lemma~1]{Hotz2016a} implies $[\aOpt{\vLambda_i}]_1[\mP_i]_{i,j}+[\aOpt{\vLambda_j}]_1[\mP_j]_{i,j}\in\interior(\sH_k)$ and, as this is a summand in $[\fPsi(\aOpt{\vLambda_\sV},\aOpt{\vMu})]_{i,j}$, that $[\fPsi(\aOpt{\vLambda_\sV},\aOpt{\vMu})]_{i,j}\in\interior(\sH_k)$. Therewith, \cite[Corollary~4]{Hotz2016a} yields $[\fPsi(\aOpt{\vLambda_\sV},\aOpt{\vMu})]_{i,j}\neq 0$ and \eqref{eqn:akdef} implies $(\aOpt{\vLambda_\sV},\aOpt{\vMu}) \notin \sA_k$.

\section{Implications of Cost Functions on LMPs}
	\label{apx:lmpimplications}

In order to support the interpretation of Theorem~\ref{thm:exactness:pospriceacb} and the intuition behind pathological price profiles, the implications of generation cost functions on LMPs shall be highlighted. To this end, assume the generator at bus~$n$ exhibits the marginal cost $\xi>0$ and features a rectangular P-Q capability region. Therefore, its valid operating points are
\begin{align}
	\sG_n = \Big\{ \vg_n\in\nR^2 :
		\aGen{\aLB{P_n}}&\leq[\vg_n]_1\leq\aGen{\aUB{P_n}}\,,\notag\\[-0.5ex]
		\aGen{\aLB{Q_n}}&\leq[\vg_n]_2\leq\aGen{\aUB{Q_n}}\,\Big\}
	\label{eq:cost:gn}
\end{align}
and its cost function $\fC_n:\sG_n\rightarrow\nR$ reads
\begin{equation}
	\fC_n(\vg_n) = \xi\ve_1\tran\vg_n + \chi
	\label{eq:cost:cn}
\end{equation}
with $\chi\in\nR$, cf. Fig.~\ref{fig:cost:3d}. Consider the OPF problem~\eqref{eqn:opf:nonconvex:full} and its dual~\eqref{eqn:opf:nonconvex:dual}. Let $\aOpt{\vg_n}\in\sG_n$ be the corresponding optimal dispatch of the generator at bus $n$ and let $\aOpt{\vLambda_n}\in\nR^2$ be the optimal dual variables of the power balance at bus~$n$. Assume that strong duality holds in~\eqref{eqn:opf:nonconvex:full}, thus $[\aOpt{\vLambda_n}]_1$ is the LMP for active power and $[\aOpt{\vLambda_n}]_2$ the LMP for reactive power.

If the generator is \emph{not} operated at any of its limits, i.e., $\aOpt{\vg_n}\in\interior(\sG_n)$, then $\fC_n$ is differentiable and
\begin{equation}
	\aOpt{\vLambda_n}
		= \nabla\fC_n(\aOpt{\vg_n})
		= \begin{bmatrix} \xi\\ 0 \end{bmatrix}.
\end{equation}
Thus, the LMP for active power is positive and the LMP for reactive power is zero. If $\aOpt{\vg_n}$ is at the boundary of the P-Q capability region, differentiability is lost and the LMPs lie in the respective subdifferential of $\fC_n$, i.e., $\aOpt{\vLambda_n}\in\subdiff\fC_n(\aOpt{\vg_n})$. The first dimension of the subdifferential reads
\begin{equation}
	[\subdiff\fC_n(\aOpt{\vg_n})]_1 =
		\begin{cases}
			[\,\xi,+\infty) & \text{if~} [\aOpt{\vg_n}]_1 = \aGen{\aUB{P_n}} \\
			(-\infty,\xi\,] & \text{if~} [\aOpt{\vg_n}]_1 = \aGen{\aLB{P_n}} \\
			\{\xi\}         & \text{otherwise}
		\end{cases}
	\label{eq:cost:subdiff:active}
\end{equation}
which follows from the supporting hyperplanes of the epigraph of $\fC_n$, cf. Fig.~\ref{fig:cost:active}. Analogously, the second dimension reads
\begin{equation}
	[\subdiff\fC_n(\aOpt{\vg_n})]_2 =
		\begin{cases}
			[\,0,+\infty) & \text{if~} [\aOpt{\vg_n}]_2 = \aGen{\aUB{Q_n}} \\
			(-\infty,0\,] & \text{if~} [\aOpt{\vg_n}]_2 = \aGen{\aLB{Q_n}} \\
			\{0\}         & \text{otherwise}
		\end{cases}
	\label{eq:cost:subdiff:reactive}
\end{equation}
see also Fig.~\ref{fig:cost:reactive}. It follows from~\eqref{eq:cost:subdiff:active} that the LMP for active power may only become zero or negative if the generator is operated at its minimum active power output. In this case, the LMP for active power typically comprises the marginal cost of another generator (augmented by the cost of transmission losses), which usually yields a positive price. Nonetheless, the LMP of active power may become negative, e.g., if additional load at bus $n$ mitigates the impact of power flow restrictions (e.g., due to congestion) and, therewith, introduces more social welfare than the generation cost it causes. For reactive power, it follows from~\eqref{eq:cost:subdiff:reactive} that its LMP may only become negative if additional inductive reactive power increases social welfare. However, this circumstance is counteracted by the inductive nature of transmission grids. In all other cases, the LMP for reactive power is zero or positive.

These considerations apply similarly to general convex and increasing generation cost functions as well as buses without any generator. For example, the latter implies ${\sG_n=\{\mZero\}}$. This is a special case of $\sG_n$ in~\eqref{eq:cost:gn}, i.e., when all bounds are zero. Correspondingly, the preceding discussion applies as well.

\section*{Acknowledgment}

M.~Hotz would like to thank Andrej Joki{\'c} of the University of Zagreb as well as Maximilian Riemensberger, Lorenz Weiland, and Niklas Winter of the Technische Universit\"at M\"unchen for valuable discussions.

\IEEEtriggeratref{13}

\bibliographystyle{IEEEtran}
\bibliography{hotz_utschick_tps2016}

\end{document}